\documentclass{article}

\usepackage{amsthm}

\usepackage[T1]{fontenc}
\usepackage[english]{babel}

\usepackage{lipsum}
\usepackage{amsfonts}
\usepackage{graphicx}
\usepackage{epstopdf}
\usepackage{algorithmic}
\ifpdf
  \DeclareGraphicsExtensions{.eps,.pdf,.png,.jpg}
\else
  \DeclareGraphicsExtensions{.eps}
\fi

\usepackage{amsmath}
\usepackage{amssymb}
\usepackage{pgfplotstable}
\usepackage{booktabs}
\usepackage{colortbl}

\usepackage{hyperref}
\usepackage{cleveref}
\newfont{\NUMBERS}{msbm9 scaled\magstep1}

\newcommand{\REAL}{ \mathbb{R}}

\usepackage{ifthen}

\newcommand{\REALP}{\REAL_{>0}}
\newcommand{\REALNN}{\REAL_{\geq 0}}

\newcommand{\Of}[1]{\left[#1\right]}
\newcommand{\Vect}[2][]
{
  \ifthenelse{\equal{#1}{}}
  {\boldsymbol{#2}}
  {#2_{#1}}
}
\newcommand{\Matr}[2][]
{
  \ifthenelse{\equal{#1}{}}
  {\boldsymbol{#2}}
  {{#2}_{#1}}
}
\newcommand{\Pinv}[1]{{#1}^{+}}

\usepackage{xargs}
\newcommand{\MatrPlus}[3][1={}, 2={}]{
  \ifthenelse{\equal{#1}{} \AND \equal{#2}{}}
  {\boldsymbol{#3}}
  {}
  \ifthenelse{\not{\equal{#1}{}} \AND \equal{#2}{}}
  {\boldsymbol{#3}_{#1}}
  {} 
  \ifthenelse{\equal{#1}{} \AND \not{\equal{#2}{}}}
  {\boldsymbol{#3}^{#2}}
  {}
  \ifthenelse{\not{\equal{#1}{}} \AND \not{\equal{#2}{}}}
  {#3_{#1,#2}}
  {}
}

\newcommand{\Kernel}[1]{\mbox{Ker}(#1)}
\newcommand{\Image}[1]{\mbox{Im}(#1)}

\newcommand{\LaplSymbol}{\Delta}
\newcommand{\Lapl}[1][]{
  \ifthenelse{\equal{#1}{}}
  {\Matr{\LaplSymbol}}
  {\Matr{\LaplSymbol}_{#1}}
}

\newcommand{\Diag}[1]{\text{{Diag}}\left(#1\right)}
\newcommand{\BlkDiag}[1]{\text{{Block Diag}}\left(#1\right)}

\newcommand{\dx}{\, dx}

\newcommand{\dt}{\, dt}

\DeclareMathOperator{\Div}{\operatorname{div}}
\DeclareMathOperator{\Grad}{\nabla}

\def\XXint#1#2#3{{\setbox0=\hbox{$#1{#2#3}{\int}$ }
\vcenter{\hbox{$#2#3$ }}\kern-.6\wd0}}

\newcommand{\Dt}[1]{\partial_{t}#1}

\newcommand{\Dim}{d}
\newcommand{\Domain}{\Omega}

\newcommand{\tstep}{k}

\newcommand{\Deltat}[1][]{
  \ifthenelse{\equal{#1}{}}
{\Delta t}
  {\Delta t_{#1}}
  }

\newcommand{\MeshPar}{h}
\newcommand{\Cell}[1][]{\ifthenelse{\equal{#1}{}}{c}{c_{#1}}}
\newcommand{\Edge}[1][]{\ifthenelse{\equal{#1}{}}{e}{e_{#1}}}

\newcommand{\SubCell}[1][]{\ifthenelse{\equal{#1}{}}{t}{t_{#1}}}
\newcommand{\Tsymb}{\mathcal{T}}
\newcommand{\Esymb}{\mathcal{E}}
\newcommand{\Triang}[1][]
{
  \ifthenelse{\equal{#1}{}}
  {\Tsymb}
  {\Tsymb_{#1}}
}
\newcommand{\Edges}[1][]
{
  \ifthenelse{\equal{#1}{}}
  {\Esymb}
  {\Esymb_{#1}}
}

\newcommand\coarse[1]{{{#1}'}}
\newcommand\fine[1]{{#1}}

\newcommand{\PotOf}[2][]
{
  \ifthenelse{\equal{#1}{}}
  {\Pot(#2)}
  {\Pot_{#1}(#2)}
}

\newcommand{\Lspacechar}{L}
\newcommand{\Lspace}[1]{\Lspacechar^{#1}}

\newcommand{\Cachar}{\mathcal{C}}

\newcommand{\Cont}[1][]{
  \ifthenelse{\equal{#1}{}}
  {\Cachar} 
  {\Cachar^{#1}}
}

\newcommand{\HolderExp}{\delta}
\newcommand{\Holder}[1][]{
  \ifthenelse{\equal{#1}{}}
  {\Cachar^{\HolderExp}}
  {\Cachar^{#1}}  
}

\newcommand{\Lip}[1]
{
  \ifthenelse{\equal{#1}{}}
  {\text{Lip}}
  {\text{Lip}_{#1}}        
}

\newcommand{\Sob}[2][]{
  \ifthenelse{\equal{#2}{}}
  {H^{#1}}
  {W^{#1,#2}}  
}

\newcommand{\DHolder}[2][]
{
  \ifthenelse{\equal{#2}{}}
  {\Cachar^{#1,\HolderExp}}   
  {\Cachar^{#1,#2}}  
}

\newcommand{\VspaceSymb}{\mathcal{V}}
\newcommand{\Vspace}[1][]{
  \ifthenelse{\equal{#1}{}}
  {\VspaceSymb_{\MeshPar}}
  {\VspaceSymb_{\MeshPar}(#1)}
}

\newcommand{\VbaseSymb}{\varphi}
\newcommand{\Vbase}[1][]{
  \ifthenelse{\equal{#1}{}}
  {\VbaseSymb}
  {\VbaseSymb_{#1}}
}
\newcommand{\WspaceSymb}{\mathcal{W}}
\newcommand{\Wspace}[1][]{
  \ifthenelse{\equal{#1}{}}
  {\WspaceSymb_{\MeshPar}}
  {\WspaceSymb_{\MeshPar}(#1)}
}

\newcommand{\WbaseSymb}{\psi}
\newcommand{\Wbase}[1][]{
  \ifthenelse{\equal{#1}{}}
  {\WbaseSymb}
  {\WbaseSymb_{#1}}
}
\newcommand{\Tdens}{\rho}
\newcommand{\Rtdens}{\tilde{\rho}}

\newcommand{\Pot}{\phi}

\newcommand{\Vel}{v}

\newcommand{\Recon}[1]{\Matr{R}_{#1}}
\newcommand{\DRecon}[1]{\Matr{R}_{#1}}

\newcommand{\Mgrad}[1]{\Matr{\nabla}_{#1}}
\newcommand{\Mdiv}[1]{\Matr{\Div}_{#1}}

\newcommand{\Source}{\Tdens^{\mbox{in}}}
\newcommand{\Sink}{\Tdens^{\mbox{f}}}

\newcommand{\BallSymb}{B}
\newcommand{\Ball}[2][]
{
  \ifthenelse{\equal{#2}{}}
  {\BallSymb_{#1}}
  {\BallSymb(#2,#1)}
}
\newcommand{\Avgint}[3][]
{
  \ifthenelse{\equal{#1}{}}
  {({#3})_{#2}}
  {(#3)_{#1,#2}}
}

\newcommand{\Norm}[2][]
{
  \ifthenelse{\equal{#1}{}}
  {{\left\lVert#2\right\rVert}}
  {{\left\lVert#2\right\rVert}_{#1}}
}

\newcommand{\Fnewton}{F}

\newcommand{\Schur}{Schur}

\newcommand{\Amat}{\mathcal{A}}
\newcommand{\Smat}{\mathcal{S}}
\newcommand{\Ablock}{\Matr{\Amat}}
\newcommand{\Bblock}{\Matr{\mathcal{B}}}
\newcommand{\BB}{\Matr{B}\Matr{B}}
\newcommand{\Cblock}{\Matr{\mathcal{C}}}

\newcommand{\Sblock}{\Matr{\mathcal{S}}}

\newcommand{\Eblock}{\Matr{\mathcal{E}}}

\newcommand{\AblockTdens}{\Matr{\mathcal{\tilde{A}}}}
\newcommand{\tildeBblock}{\Matr{\mathcal{\tilde{B}}}}

\newcommand{\Relax}{\mu}

\newcommand{\Diff}{E}

\usepackage[mathscr]{euscript}

\newcommand{\TdensToPot}[1]{\Matr{J}_{#1}}
\newcommand{\TdensToPotBlock}[1]{\Matr{\mathcal{J}}_{#1}}
\newcommand{\Identity}[1]{\Matr{I}_{#1}}

\newcommand{\IcellK}{i}
\newcommand{\IcellL}{j}
\newcommand{\Ncell}{N_{\Triang}}
\newcommand{\cncell}{N_{\coarse{\Triang}}}
\newcommand{\Iedge}{k}
\newcommand{\Nedge}{N_{\Edges}}
\newcommand{\cnedge}{N_{\coarse{\Edges}}}
\newcommand{\Ntdens}{\cncell}
\newcommand{\Npot}{\Ncell}
\newcommand{\Ntimestep}{K}
\newcommand{\Ntimep}{\Ntimestep+1}
\newcommand{\Ntime}{\Ntimestep}
\newcommand{\Slack}{s}

\newcommand{\Ortho}{|w|}

\newcommand{\Saddle}{\mathscr{J}}

\newcommand{\ZeroMat}{\ }

\newcommand{\SchurDual}{\Matr{\mathcal{S}}}

\newcommand{\InnerTol}{\varepsilon_{\text{in}}}
\newcommand{\OuterTol}{\varepsilon_{\text{out}}}

\newcommand{\TimeDiff}{\mathcal{D}_{t}}

\newcommand{\Gradblock}{\mathcal{D}_{x}}
\newcommand{\Divblock}{\mathcal{D}iv_{x}}
\newcommand{\MassTdens}{{\coarse{\Matr{M}}}}
\newcommand{\MassPot}{\Matr{M}}
\newcommand{\MassTdensBlock}{{\Matr{\coarse{\mathcal{M}}}}}
\newcommand{\MassPotBlock}{\Matr{\mathcal{M}}}
\newcommand{\TimeAvgBlock}{\Matr{\mathcal{H}}}

\newcommand{\Np}{n}
\newcommand{\Nr}{m}

\newcommand{\Approx}[1]{\hat{#1}}

\usepackage{multirow,makecell}
\newcolumntype{R}[1]{>{\raggedleft\arraybackslash }b{#1}}
\newcolumntype{L}[1]{>{\raggedright\arraybackslash }b{#1}}
\newcolumntype{C}[1]{>{\centering\arraybackslash }b{#1}}

\newcommand{\Inv}[1]{(#1)^{-1}}

\newcommand{\Block}[2]{{#2}^{#1}}
\newcommand{\Gpot}{\mathcal{G}}

\newcommand{\Shift}{\lambda}
\newcommand{\New}[1]{\bar{#1}}
\newcommand{\TPot}{\New{\Pot}}
\newcommand{\AreaTdens}{|\coarse{\Cell}|}
\newcommand{\AreaPot}{|\Cell|}

\newcommand{\ProjTdens}{\Matr{\mathcal{P}}}

\newcommand{\BlockMassTdens}{\coarse{\mathcal{M}}}
\newcommand{\DPot}{\delta \Pot}
\newcommand{\DTPot}{\delta \TPot}
\newcommand{\DTdens}{\delta \Tdens}
\newcommand{\DSlack}{\delta \Slack}
\newcommand{\DShift}{\delta \Shift}

\newcommand{\IP}{IP}

\newcommand{\AvgOuter}{\textbf{Outer/Lin.sys.}}
\newcommand{\AvgCPU}{\textbf{CPU/Lin.sys.}}
\newcommand{\AvgInner}{\textbf{Inner/Outer}}

\newcommand{\Sol}{u}

\newcommand{\RhsPot}{c}
\newcommand{\RhsTdens}{d}

\newtheorem{Prop}{Proposition}

\crefformat{prob}{Problem~(#2#1#3)}
\Crefformat{prob}{Problem~(#2#1#3)}
\crefmultiformat{prob}{Problems~(#2#1#3)}{ and~(#2#1#3)}{, (#2#1#3)}{ and~(#2#1#3)}
\Crefmultiformat{prob}{Problems~(#2#1#3)}{ and~(#2#1#3)}{, (#2#1#3)}{ and~(#2#1#3)}

\usepackage{epigraph}

\newtheorem{remark}{Remark}

\usepackage{PRIMEarxiv}

\usepackage[utf8]{inputenc} 



\title{Efficient preconditioners for solving dynamical optimal transport via interior point methods
}
\author{
  Enrico Facca\\
  University of Bergen\\
  Bergen, Norway \\
  \texttt{enrico.facca@uib.no} \\
  \And
  Gabriele Todeschi\\
  Univ. Grenoble Alpes, ISTerre\\
  F-38058 Grenoble, France\\
  \texttt{gabriele.todeschi@univ-grenoble-alpes.fr}\\
  \And
  Andrea Natale\\
  Univ. Lille, Inria, CNRS, UMR 8524\\
  Laboratoire Paul Painlev\'e\\
  Lille, France \\
  \texttt{andrea.natalea@inria.fr} \\
  \And
  Michele Benzi\\
  Scuola Normale Superiore\\
  Piazza dei Cavalieri, 7, 56126\\
  Pisa, Italy\\
  \texttt{michele.benzi@sns.it}\\
}

\begin{document}
\maketitle

\begin{abstract}
In this paper we address the numerical solution of the quadratic optimal transport problem in its dynamical form, the so-called Benamou-Brenier formulation. When solved using interior point methods, the main computational bottleneck is the solution of large saddle point linear systems arising from the associated Newton-Raphson scheme.
The main purpose of this paper is to design efficient preconditioners to solve these linear systems via iterative methods.  Among the proposed preconditioners, we introduce one based on the partial commutation of the operators that compose the dual Schur complement of these saddle point linear systems, which we refer as $\boldsymbol{B}\boldsymbol{B}$-preconditioner.
A series of numerical tests show that the $\boldsymbol{B}\boldsymbol{B}$-preconditioner is the most efficient among those presented, despite a performance deterioration in the last steps of the interior point method. It is in fact the only one having a CPU-time that scales only slightly worse than linearly with respect to the number of unknowns used to discretize the problem.

\end{abstract}

\newcommand{\sep}{\and}

\section{Introduction}
\label{sec:introduction}
Optimal transport deals with the problem of finding the optimal way to reallocate one nonnegative density into another by minimizing the total cost of displacement in space.
In recent years, numerous contributions have been made to the study of this problem, both on the theoretical and computational level. We suggest, for example, the monographs \cite{villani2003topics,santambrogio2015optimal,ambrosio2021lectures,peyre2019computational} for a detailed presentation of the subject. Due to these advances, optimal transport is nowadays an established tool for many applications including, for example, the analysis of partial differential equations (PDE)~\cite{ambrosio2005gradient},  physical modeling~\cite{santambrogio2015optimal}, data science and machine learning~\cite{peyre2019computational}, economics~\cite{galichon2016optimal}, and inverse problems~\cite{metivier2022review}.

When the cost of displacement per unit mass is given by the square of the Euclidean distance, the problem can be recast dynamically, as shown by Benamou and Brenier~\cite{benamou2000computational}. Consider a compact and convex domain $\Domain\subset \REAL^{\Dim}$ and two nonnegative densities $\Source$ and $\Sink$ in $\Lspace{1}(\Domain)$, with $\int_{\Domain}\Source\dx=\int_{\Domain}\Sink\dx$. To transport the former to the latter, we aim to find a time-dependent density $\Tdens:[0,1]\times \Domain \to \REALNN$ and a velocity field $\Vel:[0,1]\times \Domain \to \REAL^{\Dim}$ that solve the following minimization problem:
\begin{gather}
\label{prob:min-kinetic}
\min_{\Tdens,\Vel}
\int_{0}^{1}\int_{\Domain}  \dfrac{\Tdens |\Vel|^2}{2}\dt \dx
\, : \quad
\begin{cases}
  \Dt{\Tdens}+\Div(\Tdens\Vel)=0 &\text{in $[0,1]\times \Domain$} \\
  \Tdens\Vel \cdot \hat{n}=0  &\text{on $[0,1]\times \partial \Domain$}\\
  \Tdens(0,\cdot)=\Source
  \, , \, \Tdens(1,\cdot)=\Sink
\end{cases}
\,.
\end{gather}
The total kinetic energy represents the cost of displacement.
Thanks to a change of variables $(\Tdens,\Vel)\to(\Tdens,m=\Tdens\Vel)$, \eqref{prob:min-kinetic} can be rewritten as a convex optimization problem.
From the optimality conditions, one can deduce that the optimal velocity field is the gradient of a potential $\Pot:[0,1]\times \Domain \to \REAL$. The potential $\Pot$ and the density $\Tdens$ are given as the solution of the following system of PDEs:
\begin{subequations}
	\label{eq:l2}
	\begin{align}
		\label{eq:l2-pde}
		&-\Dt{\Tdens}-\Div \left(\Tdens \Grad \Pot\right) = 0, \\
		\label{eq:l2-hj}
		&\Dt{\Pot} + \frac{\|\Grad \Pot\|^2}{2} + s = 0, \\
		\label{eq:l2-ortho}
		&\Tdens \geq 0, \, s\geq 0, \, \Tdens s = 0,
	\end{align}
\end{subequations}
with boundary conditions $\Tdens(0,\cdot)=\Source, \, \Tdens(1,\cdot)=\Sink, \, \Tdens \Grad \Pot \cdot \hat{n} = 0 \text{ on } [0,1]\times\partial \Domain$. The auxiliary variable $s:[0,1]\times \Domain \to \REAL_{\geq 0}$ is related to the positivity constraint on $\Tdens$.

The Benamou-Brenier formulation offers several advantages. In addition to the total displacement cost, it also directly provides how to continuously reallocate the mass. This also constitutes a natural way to define interpolations between densities. Furthermore, it draws a clear link between optimal transport and continuum mechanics, and it is naturally suited for Eulerian discretizations. Finally, it can be easily generalized to other problems by penalizing/constraining the evolution $\Tdens$ (such as variational mean field games and planning problems \cite{achdou2012mean}, or unbalanced optimal transport \cite{chizat2018interpolating,baradat2021regularized}).
On the other hand, the numerical solution of~\eqref{prob:min-kinetic}, or its system of optimality conditions~\eqref{eq:l2}, poses significant challenges. Although it is a convex optimization problem, it is nonlinear and, in general, non-smooth for vanishing densities. Moreover, it is a time-space boundary value problem and there is a positivity constraint on $\Tdens$ to take into account.

While different strategies have been proposed to discretize the Benamou-Brenier formulation, most of these rely on staggered time discretization and a space discretization which may be based on finite differences \cite{benamou2000computational,Papadakis2014}, finite volumes \cite{gladbach2020scaling,natale2021computation,lavenant2021unconditional,haber2015multilevel} or finite elements \cite{lavenant2018dynamical,natale2022mixed,lavenant2021unconditional}. Here, we focus on the framework considered in~\cite{natale2021computation}, where one uses finite volumes in space with a two-level discretization of the domain $\Domain$, in order to discretize the density $\Tdens$ and the potential $\Pot$ separately, which alleviates some checkerboard instabilities that may appear when the same grid is used to discretize both variables (a strategy first used in~\cite{Facca-et-al:2018,Facca-et-al-numeric:2020} when dealing with optimal transport with unitary displacement cost given by the Euclidean distance). This choice is not restrictive since such a framework constitutes a generalization of the finite volume scheme studied in \cite{gladbach2020scaling, lavenant2021unconditional} (in which a single grid is used for density and potential), and it also contains as special cases the finite difference scheme in \cite{Papadakis2014} (when a single Cartesian grid is used in space), or the discrete transport models on networks studied in \cite{erbar2020computation} (by an appropriate reinterpretation of the space tessellation).

From a computational point of view, the numerical solution of the resulting discrete optimization problem is generally tackled by iterative first-order primal-dual optimization schemes (see, e.g., \cite{benamou2000computational,Papadakis2014}).
The computational cost per iteration of these methods is typically low but the total number of iterations strongly depends on the data and the tuning of the parameters, generally  increasing with the problem size. In this paper we will consider instead a second-order approach, an Interior Point (IP) method proposed in \cite{natale2021computation}, where the optimization problem~\cref{prob:min-kinetic} is relaxed by adding a logarithmic barrier function (we refer the reader to~\cite{wright1997primal,Boyd2004convopt,gondzio2012interior} for a broad introduction of IP methods). This perturbation provides smoothness by enforcing the strict positivity of the density and uniqueness of the solution. The problem can then be effectively solved by solving the relaxed optimality conditions via the Newton method, and the original unperturbed solution is retrieved by repeating this procedure while reducing the relaxation term to zero. Numerical experiments in~\cite{natale2021computation} show that the total number of Newton iterations remains practically constant, independently of the number of unknowns.

The most demanding task in \IP{} methods is the solution of a sequence of saddle point linear systems in the following form
\begin{equation*}
  \begin{pmatrix}
    \Ablock
    &
    \Bblock^T
    \\
    \Bblock
    &
    - \Cblock
  \end{pmatrix}
  \begin{pmatrix}
    \Vect{x}
    \\
    \Vect{y}
  \end{pmatrix}
  =
  \begin{pmatrix}
    \Vect{f}
    \\
    \Vect{g}
  \end{pmatrix} ,
\end{equation*}
generated by the Newton method. In our problem all matrices $\Ablock$, $\Bblock$, and $\Cblock$ change at each Newton iteration but they are sparse, hence iterative solvers are the natural candidate for solving these linear systems. Iterative methods were not considered in~\cite{natale2021computation}, and thus the goal of this paper is to determine whether one can devise a preconditioning strategy such that, ideally, the total CPU time scales linearly with the number of unknowns.

We consider three preconditioners: a preconditioner based on the approximation of the inverse of the primal Schur complement $\Sblock_p=\Ablock + \Bblock \Cblock^{-1} \Bblock^T$, the SIMPLE preconditioner~\cite{Patankar72}, and one preconditioner inspired by the ideas in~\cite{Elman2006}, where we approximate the inverse of the dual Schur complement $\Sblock=-\Cblock-\Bblock \Ablock^{-1} \Bblock^T$ looking at the differential nature of the operators that compose $\Sblock$ (note that, formally, we cannot write ${\Ablock}^{-1}$ since in our problem the matrix $\Ablock$ will be singular). A series of numerical experiments suggest that the latter, named $\BB$-preconditioner, is the most efficient among those presented in this paper, despite some loss of robustness in the last \IP{} iterations.

Let us stress that, while we conducted our experiments using the scheme described in~\cite{natale2021computation}, we do expect that the $\BB$-preconditioner will provide good results when adopting different discretization schemes even outside the framework considered here, e.g. when finite elements are used. This is because it is entirely derived from the differential operators involved in the continuous problem.

\subsection{Related works}
The literature on solving~\cref{eq:l2} using second-order and/or Newton-based methods is relatively scarce. The most closely related works are~\cite{Achdou2012,haber2015multilevel}. In the first paper, the authors study linear algebra approaches involved in the solution of mean-field games systems via the Newton method. Problem \eqref{prob:min-kinetic} may in fact be seen as the vanishing viscosity limit of a particular mean-field game. However, their linear algebra approaches do not fit into the staggered temporal grids, which is generally required to discretize the Benamou-Brenier formulation. Moreover, they lose efficiency with vanishing viscosity. In the second paper, the authors solved precisely the Benamou-Brenier problem using an inexact Newton method combined with finite volumes on Cartesian grids. However, the discretization scheme is not designed to handle densities $\Source$ and $\Sink$ with compact support. The solution of the sequence of saddle-point linear systems associated to the Newton method was performed using a preconditioning strategy that resembles the approach based on the approximation of the primal Schur complement, dropping some of the mixed time-space terms. However, this approach did not give satisfactory results in our experiments (more details in~\cref{sec:primal}).

\subsection{Paper structure}
In~\cref{sec:discretization} we summarize the discretization method proposed in~\cite{natale2021computation}, the \IP{} approach used to solve the optimization problem, and the nonlinear problems associated with it. Then, in~\cref{sec:Newton}, we present the linear system studied in this paper. Finally, in \cref{sec:preconditioning} we present the preconditioners described in this paper, together with some numerical experiments where we solve a specific test case for different time and space refinements. In particular, in~\cref{sec:resume} we compare the CPU time required by the preconditioners in solving different test cases.

\section{Discrete problem and Interior Point method}
\label{sec:discretization}
In this section we present the discrete counterpart  of system~\cref{eq:l2} considered in~\cite{natale2021computation}. This combines the use of finite volumes with staggered temporal grids. We further present the \IP{} strategy adopted for solving the problem. 

\subsection{Spatial discretization}
In~\cite{natale2021computation}, the authors considered two different discretizations of the domain $\Domain$, which is assumed to be polygonal, both admissible for Two Point Flux Approximation (TPFA) finite volumes, according to \cite[Definition 9.1]{eymard2000finite}.
At each time $t\in [0,1]$, the variable $\Tdens(t)$ is discretized on a Delaunay triangulation, with the further hypothesis that only acute angles appear. The variable $\Pot(t)$ is discretized on a finer grid obtained from the previous one by dividing each triangular cell into three quadrilateral cells, joining the edges midpoints to the triangle's circumcenter~\cite[Figure 1]{natale2021computation}. Both meshes consist of two sets $\left(\fine{\Triang},\fine{\Edges}\right)$, the set of cells $\fine{\Cell}$ and edges $\fine{\Edge}$, respectively. To distinguish the coarser mesh from the finer one, we denote the former by $\left(\coarse{\Triang},\coarse{\Edges}\right)$, with the same notation for all its elements.
Due to the no-flux boundary condition, boundary edges are not relevant to the discrete model. We will then consider, by convention, the sets $\Edges$ and $\coarse{\Edges}$ without boundary edges. Let us denote by $\cncell$ and $\cnedge$ the total number of cells and (internal) edges of the coarser mesh. The total number of cells and edges of the finer mesh are then $\fine{\Ncell}=3\cncell$ and $\fine{\Nedge}=2\cnedge+3\cncell$.

We define
\begin{gather*}
	\MassPot:=\Diag{\Vect{\AreaPot}}\ , \ \Vect{\AreaPot}:=(|\Cell[\IcellK]|)_{{\IcellK=1}}^{\Ncell},
	\\
	\MassTdens:=\Diag{\Vect{\AreaTdens}} \ , \  \Vect{\AreaTdens}:=(|\coarse{\Cell[\IcellK]}|)_{{\IcellK=1}}^{\cncell},
\end{gather*}
where $|\Cell[\IcellK]|$ denotes the area of the triangle $\Cell[\IcellK]$. Since the discrete model involves two different spatial discretizations, we also introduce the injection operator $\Matr{\TdensToPot{}}\in\REAL^{\Ncell,\cncell}$,
\begin{equation}
  \label{eq:tdens2pot}
\Matr{\TdensToPot{}}[i,j] = 
\left\{
\begin{aligned}
&1 &&\mbox{ if } \Vect[i]{\fine{\Cell}} \subset \Vect[j]{\coarse{\Cell}}\,,\\
&0 &&\mbox{ else}\,.
\end{aligned}
\right.
\end{equation}
The operator $\Matr{\TdensToPot{}}$ injects piecewise constant functions from the coarse grid to the fine one. The operator $(\Matr{\MassTdens})^{-1}\TdensToPot{}^T\Matr{\MassPot}$ is an averaging operator that acts in the opposite direction. Both operators are the identity when the same grid is used for $\Pot$ and $\Tdens$.

In order to define the discrete differential operators associated to the finite volume discretization, we fix an arbitrary orientation on the set of edges $\Edges$ and define the matrix $\Matr{\Diff}\in \REAL^{\fine{\Ncell},\fine{\Nedge}}$ given by
\begin{equation*}
\Matr{\Diff}[\IcellK,\Iedge]=
\left\{
\begin{aligned}
&1 &&\mbox{ if cell } \Cell[\IcellK] \mbox{ is the left side of edge } \Edge[\Iedge]\,,\\
-&1 &&\mbox{ if cell } \Cell[\IcellK] \mbox{ is the right side of edge } \Edge[\Iedge]\,.
\end{aligned}
\right.
\end{equation*}
The discrete gradient and divergence, $\Mgrad{} \in \REAL^{\fine{\Nedge},\fine{\Ncell}}$ and $\Mdiv{} \in \REAL^{\Ncell,\Nedge}$, are given by
\begin{gather*}
\Mgrad{}=\Mgrad{\Edges,\Triang}=\Diag{\Vect{\Ortho}}^{-1}\Matr{\Diff}^T,
\\
\Mdiv{}=\Mdiv{\Triang,\Edges}
=-\Mgrad{}^T\Diag{\Vect{\Ortho}\odot |\Vect{\Edge}|}
=-\Matr{\Diff}\Diag{|\Vect{\Edge}|},
\end{gather*}
where $\Vect{\Ortho}\in \REAL^{\Nedge}$ and $\Vect{|\Vect{\Edge}|}\in \REAL^{\Nedge}$ are the vectors of distances between the circumcenters of adjacent cells and the lengths of the edges, respectively (we used the symbol $\odot$ to denote the pointwise multiplication of vectors).

The TPFA finite volume discretization and the two-level grids require the introduction of two additional operators to match the dimension between the different spaces.  The first is a reconstruction operator $\Recon{\Edges}:\REAL^{\cncell}   \rightarrow \REAL^{\Nedge}$, mapping a positive density $\Vect{\rho}$ defined on the cells of the coarse grid to a variable defined on the edges of the finer one. In~\cite{natale2021computation} this reconstruction first lifts the density $\Vect{\rho}$ into the finer space via the operator $\Matr{\TdensToPot{}}$ and then averages the values of adjacent cells.
The authors considered two types of averages, a weighted arithmetic mean and a weighted harmonic mean (the latter is not considered in this paper to avoid complicating the exposition). Using the former, the linear reconstruction operator explicitly writes
\begin{equation*}
	(\Recon{\Edges} \Vect{\Tdens})_{\Edge}:=
	\lambda_{\Edge} (\Matr{\TdensToPot{}} \Vect{\Tdens})_{\Cell[\IcellK]}+ (1-\lambda_{\Edge}) (\Matr{\TdensToPot{}} \Vect{\Tdens})_{\Cell[\IcellL]}
\end{equation*}
for the two cells $\Cell[\IcellK],\Cell[\IcellL]$ sharing the edge $\Edge$, where $\lambda_{\Edge}$ is a weight that depends on the mesh geometry. Moreover, they introduced a further operator $\DRecon{\Triang}: \REAL^{\Nedge} \rightarrow \REAL^{\cncell}$ that maps the variables defined at the edges of the finer grid to the variables defined at the cells of the coarser one. To preserve the variational structure of the discrete problem, it is defined as
\begin{align*}
	\DRecon{\Triang}:=
	\left(  \Recon{\Edges} \right)^T \Diag{\Vect{\Ortho} \odot |\Vect{\Edge}|}.
\end{align*}

\subsection{Temporal discretization}
In \cite{natale2021computation}, the temporal discretization is based on staggered temporal grids. The time interval $[0,1]$ is divided into
$\Ntime+1$ sub-intervals $[t^{\tstep},t^{\tstep+1}]$ of equal length $\Deltat=1/(\Ntimep)$, for $\tstep=0,\ldots,\Ntime\geq1$, with $t^0=0$ and $t^{\Ntimep}=1$. The variables $\Tdens$ and $\Slack$ are defined at time $t^{\tstep}$ for $\tstep=1,\ldots,\Ntime$, while $\Pot$ is discretized at each instant $(t^{\tstep}+t^{\tstep+1})/2$ for $\tstep=0,\ldots,\Ntime$.

Combining this temporal and spatial discretization, the discrete counterpart of the potential $\Pot$, the density $\Tdens$, and the slack variable $\Slack$ in equation~\cref{eq:l2} are the vectors $\Vect{\Pot}\in \REAL^{\Np}, \Vect{\Tdens} \in \REALNN^{\Nr}$, and $\Vect{\Slack} \in \REALNN^{\Nr}$ with $\Np=\Npot  (\Ntimep)$ and $\Nr=\Ntdens \Ntime$ given by
\begin{alignat*}{2}
\Vect{\Pot}&=\left(\Vect{\Pot}^1;\ldots;\Vect{\Pot}^{\Ntimep}\right),
&
\ \Vect{\Pot}^\tstep
\in \REAL^{\Npot}, \\
\Vect{\Tdens}&=\left(\Vect{\Tdens}^1;\ldots;\Vect{\Tdens}^{\Ntime}\right),
&\ \Vect{\Tdens}^\tstep
\in \REALNN^{\Ntdens}, \\
\Vect{\Slack}&=\left(\Vect{\Slack}^1;\ldots;\Vect{\Slack}^{\Ntime}\right),
&\ \Vect{\Slack}^\tstep
\in \REALNN^{\Ntdens},
\end{alignat*}
where we use the symbol $;$ to denote the concatenation of vectors. Moreover, we will denote by $\Block{\tstep}{\Vect{x}}$ the $\tstep$-slice of a concatenated vector $\Vect{x}$, which corresponds to its $\tstep$-th time portion. 
The discrete counterparts of the initial and final densities are the two vectors $\Vect{\Tdens}^0$ and $\Vect{\Tdens}^{\Ntimep}$ in $\REAL^{\cncell}$ given by
\begin{equation*}
  \Vect[\IcellK]{\Tdens}^{0}=
  \frac{1}{| \coarse{\Cell[\IcellK]}|}
  \int_{\coarse{\Cell[\IcellK]}}\Source\dx
  \,, \ 
  \Vect[\IcellK]{\Tdens}^{\Ntimep}=
  \frac{1}{| \coarse{\Cell[\IcellK]}|} 
  \int_{\coarse{\Cell[\IcellK]}}\Sink\dx
  ,\quad \IcellK=1,\ldots,\cncell\,.
\end{equation*}

\subsection{Discrete nonlinear system of equations}
\label{sec:optimization}
With the aforementioned discretization, the discrete counterpart of the nonlinear system~\cref{eq:l2} is given by
\begin{equation}
  \label{eq:l2-discrete-mu0}
  \begin{aligned}
      \Fnewton_{\Vect{\Pot}}(\Vect{\Pot},\Vect{\Tdens})&:=
      \left(
      \Fnewton^{1}_{\Vect{\Pot}};
      \ldots;
      \Fnewton^{\Ntimep}_{\Vect{\Pot}}
      \right)& =&\Vect{0}\in \REAL^{\Npot (\Ntimep)}\ ,\\
      \Fnewton_{\Vect{\Tdens}}(\Vect{\Pot},\Vect{\Tdens},\Vect{\Slack})&:=
      \left(
      \Fnewton^{1}_{\Vect{\Tdens}};
      \ldots;
      \Fnewton^{\Ntime}_{\Vect{\Tdens}}
      \right)& =&\Vect{0}\in \REAL^{\Ntdens \Ntime}\ ,\\
      \Fnewton_{\Vect{\Slack}}(\Vect{\Tdens},\Vect{\Slack})&:=
      \left(
      \Fnewton^{1}_{\Vect{\Slack}};
      \ldots;
      \Fnewton^{\Ntime}_{\Vect{\Slack}}
      \right)& =&\Vect{0}\in \REAL^{\Ntdens \Ntime}\ ,
  \end{aligned}
\end{equation}
where for $\tstep=1,\ldots,\Ntimep$ the functions $\Fnewton^{\tstep}_{\Vect{\Pot}}$ are given by
\begin{equation}
	\textstyle
  \label{eq:conservation_discrete}
	\Fnewton^{\tstep}_{\Vect{\Pot}}(
	{\Vect{\Pot}},
	{\Vect{\Tdens}}
	)
	=
	-
        \MassPot
	\Matr{\TdensToPot{}}
	\left(
	\frac{\Vect{\Tdens}^{\tstep}-\Vect{\Tdens}^{\tstep-1}}{\Deltat}
	\right)
	-
	\Mdiv{}
	\left(
	\Recon{\Edges}\left(\frac{\Vect{\Tdens}^{\tstep}+\Vect{\Tdens}^{\tstep-1}}{2}\right)
	\odot
	\Mgrad{}
	\Vect{\Pot}^{\tstep}
	\right),
\end{equation}
while for $\tstep=1,\ldots,\Ntime$ the functions $\Fnewton^{\tstep}_{\Vect{\Tdens}}$ and $\Fnewton^{\tstep}_{\Vect{\Slack}}$ are given by
\begin{gather}
	\textstyle
  \label{eq:hj_discrete}
  \Fnewton^{\tstep}_{\Vect{\Tdens}}(
	  {\Vect{\Pot}},
	  {\Vect{\Tdens}},
	  {\Vect{\Slack}}
	  )
	  =
	  \MassTdens
	  \Matr{\TdensToPot{}}^{T} 
	  \left(
	  \frac{
	    \Vect{\Pot}^{\tstep+1}-\Vect{\Pot}^{\tstep}
	  }{\Deltat}
	  \right)
	  +
	  \frac{1}{4}
	  \DRecon{\Triang}
	  \left(\big(\Mgrad{} \Vect{\Pot}^{\tstep}\big)^2 +  \big(\Mgrad{} \Vect{\Pot}^{\tstep+1}\big)^2
	  \right)
	  +
	  \MassTdens \Vect{\Slack}^{\tstep},
	  \\[0.5em]
	  \label{eq:ortho_discrete}
	  \Fnewton^{\tstep}_{\Vect{\Slack}}(
		  {\Vect{\Tdens}},
		  {\Vect{\Slack}}
		  )
		  =
		  \Vect{\Tdens}^{\tstep}\odot \Vect{\Slack}^{\tstep},
\end{gather}
with the additional (component-wise) requirement $ \Vect{\Tdens},\Vect{\Slack}\geq \Vect{0}$ . System~\cref{eq:l2-discrete-mu0} is the system of optimality conditions for the discrete counterpart of the problem in~\cref{prob:min-kinetic}.

\begin{remark}
  \label{remark:one-mass}
  Thanks to the finite volume discretization, the conservative structure of the continuity equation is preserved. Then, due to the no-flux boundary conditions, which are encoded explicitly in the divergence operator, the equation $\Block{\tstep}{\Fnewton}_{\Vect{\Pot}}(\Vect{\Pot},\Vect{\Tdens})=\Vect{0}$ implies
  \[
  \Vect{1}^T\Fnewton^{\tstep}_{\Vect{\Pot}}=
  \Vect{\AreaPot}^{T} 
  \Matr{\TdensToPot{}}
  \left(
  \Vect{\Tdens}^{\tstep}-\Vect{\Tdens}^{\tstep-1}
  \right) = 
  \Vect{\AreaTdens}^{T}
  \left(
  \Vect{\Tdens}^{\tstep}-\Vect{\Tdens}^{\tstep-1}
  \right)=0 \quad \forall \tstep=1,\ldots,\Ntime\,.
  \]
  Hence, at each intermediate time step, the discrete mass $\Vect{\AreaTdens}^{T}\Vect{\Tdens}^{\tstep}$ is preserved and is equal to the mass of the discrete initial and final densities.
\end{remark}

\subsection{Interior Point method}
Due to the \IP{} strategy, the discrete optimization problem is perturbed by adding a logarithmic barrier function,
\[
-\mu\sum_{k=1}^{\Ntime} \Deltat\sum_{i=1}^{\Ntdens} \log(\Tdens^k_i)|c_i|
\]
tuned by a parameter $\mu>0$.
This turns the system of optimality conditions~\cref{eq:l2-discrete-mu0} into
\begin{equation}
	\label{eq:l2-discrete}
	\begin{aligned}
	  \Fnewton_{\Vect{\Pot}}(\Vect{\Pot},\Vect{\Tdens})&=
	  \Vect{0}\in \REAL^{\Npot (\Ntimep)}\ ,\\
	  \Fnewton_{\Vect{\Tdens}}(\Vect{\Pot},\Vect{\Tdens},\Vect{\Slack})&=
	  \Vect{0}\in \REAL^{\Ntdens \Ntime}\ ,\\
	  \Fnewton_{\Vect{\Slack}}(\Vect{\Tdens},\Vect{\Slack})
	  -\Relax\Vect{1} &= \Vect{0} \in \REAL^{\Ntdens \Ntime}\ ,
	\end{aligned}
\end{equation}
where the complementarity constraint between $\Tdens$ and $\Slack$ in~\cref{eq:ortho_discrete} is now relaxed. In this way, the two variables are forced to be strictly positive, making the problem easier to solve by using Newton methods.

The solution of the original problem is recovered in the limit $\Relax\rightarrow 0$ and the parameter $\Relax$ directly provides a bound on the sub-optimality of the discrete solution of \cref{eq:l2-discrete} for the original unperturbed problem \cite[Section 4]{natale2021computation}. From the practical point of view, it is computed  via a continuation method: for each value of $\Relax^l$ of a sequence of relaxation parameters $(\Relax^l)_{l\ge0}$, the problem is solved with a Newton method initialized with the previously computed solution. We refer to \cite{natale2021computation} for the precise set up of the algorithm. We stress only that for each value of $\Relax$, problem~\cref{eq:l2-discrete} is solved up to the relatively low tolerance $1e-6$ in order to obtain a more robust implementation. Furthermore, within each Newton cycle, we ensure that any Newton update does not introduce negative entries in the vectors $\Vect{\Tdens}$ and $\Vect{\Slack}$ with a simple line-search procedure.

\section{The linear algebra problem and the KKT system}\label{sec:Newton}
The nonlinear system of equations~\cref{eq:l2-discrete} is solved using an inexact Newton method.  Each Newton iteration requires the solution of a linear system (referred to as KKT system) in the form
\begin{align}
	\label{eq:newton-linsys}
	\begin{pmatrix}
		\Ablock
		&
		\Bblock^T
		&
		\Matr{\ZeroMat}
		\\
		\Bblock
		&
		\Matr{\ZeroMat}
		&
		\MassTdensBlock
		\\
		\Matr{\ZeroMat}
		&
		\Diag{\Vect{\Slack}}
		&
		\Diag{\Vect{\Tdens}}
	\end{pmatrix}
	\begin{pmatrix}
		\Vect{\DPot}
		\\
		\Vect{\DTdens}
		\\
		\Vect{\DSlack}
	\end{pmatrix}
	&
	=
	\begin{pmatrix}
		\Vect{f}
		\\
		\Vect{g}
		\\
		\Vect{h}
	\end{pmatrix}
	=
	-
	\begin{pmatrix}
		\Vect{\Fnewton}_{\Vect{\Pot}}
		\\
		\Vect{\Fnewton}_{\Vect{\Tdens}}
		\\
		\Vect{\Fnewton}_{\Vect{\Slack}}-\Relax\Vect{1}
	\end{pmatrix},
\end{align}
where the block matrix in equation~\cref{eq:newton-linsys} is the Jacobian matrix of $\left( \Vect{\Fnewton}_{\Vect{\Pot}}; \Vect{\Fnewton}_{\Vect{\Tdens}}; \Vect{\Fnewton}_{\Vect{\Slack}}-\Relax\Vect{1}\right)$.
We denote by $\Matr{\MassTdensBlock}$ and $\Matr{\MassPotBlock}$ the matrices given by
\begin{equation*}
  \label{eq:masses-block}
  \Matr{\MassTdensBlock}=\BlkDiag{(\MassTdens)_{\tstep=1}^{\Ntime}},
  \quad
  \Matr{\MassPotBlock}=\BlkDiag{(\MassPot)_{\tstep=1}^{\Ntimep}}.
\end{equation*}
The matrices $\Ablock, \Bblock,\Bblock^T$ in~\cref{eq:newton-linsys} are the finite-dimensional version of the following differential operators (up to a multiplication by a mass matrix)
\begin{equation*}
  \Ablock\approx-\Div(\Tdens \Grad)
  \,,\quad
  \Bblock \approx\Dt{}+\Grad \Pot \cdot \Grad
  \,,\quad
  \Bblock^T\approx -\Dt{}-\Div( \cdot \Grad \Pot)\,.
\end{equation*}

According to the discretization described in~\cref{sec:discretization}, matrix $\Ablock\in\REAL^{\Np,\Np}$ is a symmetric block diagonal matrix with $\Ntimep$ blocks. Each block is a weighted Laplacian matrix given by
\begin{align*}
  \Block{\tstep}{\Matr{A}}
  &=
  -
  \Mdiv{}
  \Diag{\Vect{\Rtdens}^{\tstep}
  }
  \Mgrad{}\in \REAL^{\Npot,\Npot},
  \quad
  \Vect{\Rtdens}^{\tstep}:=
  \Recon{\Edges}\left( \frac{\Vect{\Tdens}^{\tstep}+\Vect{\Tdens}^{\tstep-1}}{2} \right),
\end{align*}
which we can write equivalently as
$\Ablock =-\Matr{\Divblock} \Diag{
  (\Vect{\Rtdens}^1;
  \ldots; \Vect{\Rtdens}^{\tstep};
  \ldots;  \Vect{\Rtdens}^{\Ntimep}
  )
} \Matr{\Gradblock}
$ where
\begin{equation*}
  \label{eq:divblock-gradblock}
  \Matr{\Divblock}
  =
  \BlkDiag{(\Mdiv{})_{\tstep=1}^{\Ntimep}}
   , \quad 
  \Matr{\Gradblock}
  =
  \BlkDiag{(\Mgrad{})_{\tstep=1}^{\Ntimep}}.
\end{equation*}

Matrix $\Bblock$ is a block bi-diagonal matrix in $\REAL^{\Nr, \Np}$ given by
\begin{align*}
  \Bblock
  &=
  \TdensToPotBlock{}^T\Matr{\TimeDiff}\MassPotBlock{}
  +\TimeAvgBlock\Matr{\Gpot}
    \Matr{\Gradblock}
    \,,
\end{align*}
where the matrices $\Matr{\TimeDiff}\in \REAL^{\Ntime\!\Npot,(\Ntime\!+1)\Npot}$, $\Matr{\TimeAvgBlock}\in \REAL^{\Ntime\!\Ntdens,(\Ntime\!+1)\Ntdens}$, $\TdensToPotBlock{}\in \REAL^{(\Ntime\!+1)\Npot,(\Ntime\!+1)\Ntdens}$, and $\Matr{\Gpot}\in \REAL^{(\Ntimep)\Ntdens,(\Ntimep)\coarse{\Nedge}}$ are
\begin{gather}
  \label{eq:dt-it}
  \Matr{\TimeDiff}
  \!
  =
  \!
  \frac{1}{\Deltat}
 \begin{pmatrix}
    -\Identity{\Npot}
    &
    \Identity{\Npot}
    &
    \ZeroMat
    &
    \ZeroMat
    \\
    \ZeroMat
    &
    \ddots
    &
    \ddots
    &
    \ZeroMat
    \\
    \ZeroMat
    &
    \ZeroMat
    &
    -\Identity{\Npot}
    &
    \Identity{\Npot}
  \end{pmatrix}
  \, ,\quad 
  \TdensToPotBlock{}
  =
  \BlkDiag{(\TdensToPot{})_{\tstep=1}^{\Ntimep}}
  , \\
  \label{eq:timeavg-G}
  \Matr{\TimeAvgBlock}
  =
  \frac{1}{2}
  \begin{pmatrix}
    \Identity{\Ntdens}
    &
    \Identity{\Ntdens}
    &
    \ZeroMat
    &
    \ZeroMat
    \\
    \ZeroMat
    &
    \ddots
    &
    \ddots
    &
    \ZeroMat
    \\
    \ZeroMat
    &
    \ZeroMat
    &
    \Identity{\Ntdens}
    &
    \Identity{\Ntdens}
  \end{pmatrix}
  \, , \quad
  \Matr{\Gpot}
  =
  \BlkDiag{  (\Block{\tstep}{\Matr{G}})_{\tstep=1}^{\Ntimep}},
\end{gather}
with $(\Block{\tstep}{\Matr{G}})_{\tstep=1,\ldots,\Ntimep}\in\REAL^{\Ntdens,\fine{\Nedge}}$ given by
\begin{equation*}
  \Block{\tstep}{\Matr{G}}:=
  \DRecon{\Triang}
  \Diag{\Mgrad{} \Vect{\Pot}^{\tstep}}.
\end{equation*}

\begin{remark}
  \label{remark:kernel}
  The kernel of the block diagonal matrix $\Ablock$ has dimension $\Ntimep$, since each block $\Block{\tstep}{\Matr{A}}$ has the constant vectors as kernel.
  On the other hand, $\Kernel{\Ablock}$ is included in $\Kernel{\TimeAvgBlock\Matr{\Gpot}
  \Matr{\Gradblock}}$ and therefore
  \[
  \Kernel{\Ablock}\cap \Kernel{\Bblock} = \Kernel{\Ablock}\cap \Kernel{
    \TdensToPotBlock{}^T\Matr{\TimeDiff}\MassPotBlock{}} = \langle \Vect{1} \rangle \in
  \REAL^{\Np},
  \]
  thus the Jacobian matrix is singular. This singularity can be removed grounding one entry of the solution $\Vect{\DPot}$ or via regularization techniques~\cite{BochevLehoucq2005}. However, we prefer avoiding these approaches in this paper, since iterative methods work also when dealing with  singular matrices~\cite{van2003iterative,elman2014finite}.
\end{remark}

\subsection{Reduction to a saddle point linear system}
\label{sec:linear-solver}
In order to solve the linear system in~\cref{eq:newton-linsys} it is convenient to eliminate the variable $\Vect{\DSlack}$ writing
\[
\Vect{\DSlack}=\Inv{\Diag{\Vect{\Tdens}}}(\Vect{h}-\Diag{\Vect{\Slack}}{\Vect{\DTdens}})
\]
and reduce it into the following saddle point system
\begin{align}
  \label{eq:saddle-point}
  \Matr{\Saddle}
  \begin{pmatrix}
    \Vect{\DPot}
    \\
    \Vect{\DTdens}
  \end{pmatrix}
  =
  \begin{pmatrix}
    \Ablock
    &
    \Bblock^T
    \\
    \Bblock
    &
    - \Cblock
  \end{pmatrix}
  \begin{pmatrix}
    \Vect{\DPot}
    \\
    \Vect{\DTdens}
  \end{pmatrix}
  &
  =
  \begin{pmatrix}
    \Vect{f}
    \\
    \tilde{\Vect{g}}
  \end{pmatrix},
\end{align}
where the matrix $\Cblock\in \REAL^{\Nr,\Nr}$ and the vector $\tilde{\Vect{g}}\in \REAL^{\Nr}$ are given by
\begin{gather*}
	\Cblock:=
	\MassTdensBlock
	\Diag{\Vect{\Tdens}}^{-1}
	\Diag{\Vect{\Slack}} 
	\,,\,
	\tilde{\Vect{g}}:=\Vect{g}-\MassTdensBlock\Diag{\Vect{\Tdens}}^{-1}\Vect{h} \,.
\end{gather*}
The matrix $\Matr{\Saddle}$ in~\cref{eq:saddle-point} has typically higher conditioning number than the Jacobian matrix in \cref{eq:newton-linsys}, in particular when approaching the optimal solution, i.e. as~$\Vect{\Slack}\odot \Vect{\Tdens}\approx\mu \Vect{1}\to\Vect{0}$. However, designing an efficient preconditioner for the fully coupled system in~\cref{eq:newton-linsys} can be harder than for the standard saddle point linear system in~\cref{eq:saddle-point}, for which different preconditioning approaches are described in~\cite{Benzi-et-al:2005,bergamaschi2004preconditioning,morini2017comparison}.

We adopt an Inexact Newton approach, which means that we seek for a solution $(\Vect{\DPot},\Vect{\DTdens})$ such that
\begin{equation*}
  \|
  \Matr{\Saddle}
  (\Vect{\DPot};\Vect{\DTdens})
  -
  ( \Vect{f};\Vect{\tilde{g}})
  \|
  \leq \OuterTol
  \|(
  \Vect{f};
  \Vect{g};
  \Vect{h}
  )\| \,.
\end{equation*}
The linear system residual is scaled by the right-hand side of the original system to avoid over-solving issues we faced using $\|(\Vect{f}; \Vect{\tilde{g}})\|$.
In our experiments we use a fixed tolerance $\OuterTol=1e-5$.
This value may be quite small compared to the literature for an inexact Newton approach (for example \cite{bellavia1998inexact,gondzio2013convergence}), but we want to avoid more aggressive inexact strategies which could undermine the effectiveness of the nonlinear solver and focus on the linear algebra only. 

\begin{remark}
	\label{remark:zero-mean}
	From~\cref{remark:one-mass}, we know the solution $\Vect{\Tdens}$ to equation~\cref{eq:l2-discrete} must satisfy
	$\Vect{\AreaTdens}^{T} \Block{\tstep}{\Vect{\Tdens}}=\Vect{\AreaTdens}^{T}\Vect{\Tdens}^0 $
	for all
	$\tstep=1,\ldots\Ntime$. If this condition is satisfied for the starting point of the Newton scheme, then the increment $\Vect{\DTdens}$ solution to equations \eqref{eq:newton-linsys} satisfies
	\begin{equation*}
		\label{eq:zero-mean}
		\Vect{\AreaTdens}^{T} \Block{\tstep}{\Vect{\DTdens}}=0 \quad \forall \tstep=1,\ldots\Ntime\,,
	\end{equation*}
	at each Newton iteration.
	To see this, sum both sides of the first block of equations of system \eqref{eq:saddle-point}, for each time step $k$, and use the no-flux boundary conditions in the divergence operator.
	Due to our inexact Newton steps, this condition is not verified exactly, but only up to the tolerance chosen. We impose it by renormalizing $\Vect{\Tdens}$ after each Newton update.
\end{remark}

\section{Preconditioning approaches and numerical results}
\label{sec:preconditioning}
In this section we present three preconditioners to solve the linear system in~\cref{eq:saddle-point} and a series of numerical experiments used to measure their performance. All of these preconditioners involve the usage of inexact inner solvers, therefore, they may not be linear and must be applied as right preconditioners within a FGMRES cycle \cite{Saad1993}. In all cases, the inner solver we adopted is the AGgregation-based Algebraic MultiGrid (AGMG) method described in \cite{YvanNotay2010}.

We adopt three measures to evaluate the preconditioner performance for each Newton cycle associated with each \IP{} iteration:
\begin{enumerate}
\item
\label{item:avgouter}
\AvgOuter:{ the average number of outer FGMRES iterations. It is calculated as the number of outer iterations divided by the number of linear systems solved for each Newton cycle.}
\item
\label{item:avgcpu}
\AvgCPU:{ the average CPU time (measured in seconds) required for the solution of the linear systems. It is computed as the total CPU time required to solve all linear systems within a Newton cycle divided by the number of linear systems solved (all experiments were conducted single core, on a machine equipped with an Intel$^{\circledR}$ Xeon$^{\circledR}$ 2.8 GHz CPU and 128 GiB of RAM memory).}
\item
\label{item:avginner}
\AvgInner:{ the average number of inner iterations per preconditioner application. It is computed by dividing the cumulative number of AGMG iterations required to solve the linear systems involved in the preconditioner application (for each preconditioner, we will specify the linear systems to which we will refer) by the total number of outer iterations for each Newton cycle.}
\end{enumerate}

\subsection{Test cases}
\label{sec:test-cases}
The \IP{} algorithm and the linear solver strategies described in the next sections are tested on four numerical experiments adapted from~\cite[Section 5.2]{natale2021computation}. The first three are set in $\Domain=[0,1]^2$ whereas the fourth one, three dimensional, in $\Domain=[0,1]^3$. They are:
\begin{enumerate}
\item Gaussian densities $\Source$ and $\Sink$ centered at $(0.3,0.3)$ and $(0.7,0.7)$ respectively with variance $0.1$. This ensures that $\Source$ and $\Sink$ are lower bounded by $1e-3$.
\item Translation of a compactly supported smooth sinusoidal density $\Source$.
\item Compression of a compactly supported smooth sinusoidal density $\Source$. This test case has been shown experimentally in \cite{natale2021computation} to exhibit severe instabilities in the discrete solution, and this, in turn, has a negative effect on the discrete solver.
\item Same as test 1 but in three dimensions. The Gaussian densities $\Source$ and $\Sink$ are centered at $(0.3,0.3,0.3)$ and $(0.7,0.7,0.7)$.
\end{enumerate}

Test cases 1, 2 and 3 are discretized using four different meshes, taken from~\cite{FVCAV}, having $\cncell=224,896,3584,14336$ cells, where each mesh is a refinement of the other. These meshes will be denoted by $\Triang^0,\Triang^1,\Triang^2,\Triang^3$. The corresponding number of subcells is $\Ncell=672,2688,10752,43008$.  For each grid, we consider $\Deltat=1/(\Ntime+1)$ with $\Ntime+1=16,32,64,128$. The typical mesh size $\MeshPar$ of the mesh $\Triang^0$  is approximately $1/16$, so that the diagonal pairing of the space and time steps provides a uniform discretization of the time-space domain $[0,1]\times\Domain$. The total number of $\Pot$ and $\Tdens$ degrees of freedom for these sequence of problems is $\Np+\Nr\approx(1.4e4,1.1e5,9.1e5,7.3e6)$, which is the size of the linear system in~\cref{eq:saddle-point}. 
The last test case is solved on Cartesian grids with cubic cells. In this case the operator $\Matr{\TdensToPot{}}$ in\cref{eq:tdens2pot} is simply the identity. We consider 2 refinements (both in time and in space) of an initial experiment set on a $8\times 8 \times 8$ grid and using $\Deltat=1/(\Ntime+1)$ with $\Ntime+1=8$. In this case, the total number of $\Pot$ and $\Tdens$ degrees of freedom for a uniform time-space refinement is $\Np+\Nr\approx(7.7e3,1.3e5,2.0e6)$.

Convergence is achieved when the relaxation parameter $\mu$ is below $1e-6$, which corresponds to ten iterations of IP, with the choice of parameters described in~\cite{natale2021computation}. The nonlinear system~\cref{eq:l2-discrete} is solved with rather tight tolerance $1e-6$ in order to have a robust implementation for the IP solver. Each IP step requires between $3$ and $7$ inexact Newton iterations. These numbers are practically insensitive to the mesh and the time step adopted. The total number of Newton steps (which corresponds to the number of linear systems to be solved) ranges between $40$ and $50$. The preconditioning approaches used to solve these linear systems are described in the following sections.

For each preconditioner, we will present a table reporting the metrics mentioned above, while the relaxation parameter $\Relax$ is reduced for all combinations of spatial and temporal discretizations. We will show these detailed results only for the second test case, since it captures the main challenges involved in the linear algebra problems. The numerical results for the other three cases will be summarized in~\cref{sec:resume}, where we will compare the performance of the preconditioners in terms of total CPU time.

The code is written in \textsc{Matlab}$^{\circledR}$ and its source code is available at the link \href{https://github.com/gptod/OT-FV}{https://github.com/gptod/OT-FV}, where it is possible to reproduce the experiments presented in this paper. The AGMG solver is interfaced via \textsc{MEX}$^{\circledR}$.

\subsection{Preconditioner based on the primal Schur complement}
\label{sec:primal}
In this section we consider a preconditioner for the linear system in~\cref{eq:saddle-point} based on the primal Schur complement $\Matr{\Sblock}_p=\Ablock+\Bblock^T\Cblock^{-1} \Bblock$, similarly to the approach used in~\cite{FaccaBenzi:2021} for the solution of the optimal transport problem on graphs with cost equal to the shortest path distance.  It is given by
\begin{gather}
  \label{eq:sac-prec}
  \Matr{P}^{-1}
  =
  \begin{pmatrix}
    \Matr{I}& \ZeroMat
    \\
    \Matr{\Cblock}^{-1}\Bblock& \Matr{I}
  \end{pmatrix}
  \begin{pmatrix}
    \Matr{\Sblock}_p^{-1}& \ZeroMat
    \\
    \ZeroMat& -\Matr{\Cblock}^{-1}
  \end{pmatrix}
  \begin{pmatrix}
    \Matr{I}& \Bblock^{T}\Matr{\Cblock}^{-1}
    \\
    \ZeroMat& \Matr{I}
  \end{pmatrix}.
\end{gather}

The application of the preconditioner in~\cref{eq:sac-prec} requires the inversion of $\Matr{\Cblock}$, which is diagonal, and the (approximate) solution of the linear system
\begin{equation}
  \label{eq:primal-linsys}
  \Matr{\Sblock}_p \Vect{x} = \Vect{b}\,,
\end{equation}
with $\Vect{x},\Vect{b}\in \REAL^{\Np}$. Using this preconditioner is essentially equivalent to reducing the linear system in~\cref{eq:saddle-point} to the variable $\Vect{\DPot}$ only, an approach referred to in the literature as fully reduced \cite{behata:2009} or condensed \cite{DApuzzo2010}. 
The matrix $\Matr{\Sblock}_p$ is a block tridiagonal matrix. It can be written as
\begin{equation*}
  \label{eq:pde-primal}
  \begin{split}
  \Matr{\Sblock}_p=
  \overbrace{-\Matr{\Divblock} \Diag{\Vect{\Rtdens}} \Matr{\Gradblock}}^{\Ablock}
  \overbrace{-\Matr{\Divblock} \Matr{\Gpot}^T \TimeAvgBlock^T \Matr{\Cblock}^{-1}
  \TimeAvgBlock\Matr{\Gpot} \Matr{\Gradblock}}^{=:\Matr{\Sblock}_{xx}}
  +
  \overbrace{
    \MassPotBlock
    \Matr{\TimeDiff}^{T}    
    \TdensToPotBlock{}
    \Matr{\Cblock}^{-1}
    \TdensToPotBlock{}^T
    \Matr{\TimeDiff}
    \MassPotBlock
  }^{=:\Matr{\Sblock}_{tt}}
  \\
  +\underbrace{
    \MassPotBlock
    \Matr{\TimeDiff}^{T}    
    \TdensToPotBlock{}
    \Matr{\Cblock}^{-1}
    \TimeAvgBlock \Matr{\Gpot} \Matr{\Gradblock}
    +
    (
    \MassPotBlock
    \Matr{\TimeDiff}^{T}    
    \TdensToPotBlock{}
    \Matr{\Cblock}^{-1}
    \TimeAvgBlock \Matr{\Gpot} \Matr{\Gradblock}
    )^T}_{
    =:\Matr{\Sblock}_{tx}+\Matr{\Sblock}_{tx}^T},
  \end{split}
\end{equation*}
and it may be seen as the discretization of an elliptic time-space operator. In fact, it is  the sum of a $\Tdens$-weighted spatial Laplacian matrix $\Ablock$, an anisotropic spatial Laplacian $\Matr{\Sblock}_{xx}\approx -\Div(\Tdens/\Slack\Grad \Pot \otimes \Grad \Pot \Grad)$, a weighted temporal Laplacian with Neumann boundary condition $\Matr{\Sblock}_{tt}$, and a time-space operator $\Matr{\Sblock}_{tx}+\Matr{\Sblock}_{tx}^T$. We adopt the AGMG solver to solve the linear system in~\cref{eq:primal-linsys}, with relative accuracy $\InnerTol=1e-1>\OuterTol=1e-5$. The value of the inner tolerance $\InnerTol$ has been determined experimentally to distribute the work load between the inner and outer solver in a balanced way.

\addtolength{\tabcolsep}{-3.5pt}
\begin{table}
   \begin{center}
    {\footnotesize
      

\begin {tabular}{R{0.8cm }|R{0.45cm }R{0.45cm }R{0.45cm }R{0.52cm}|R{0.8cm}R{0.82cm}R{0.82cm}R{0.7cm}|R{0.7cm }R{0.7cm }R{0.7cm }R{0.7cm }}
\toprule
 
    $\Ntime\!+\!1$ 
    & 16 & 32 & 64 & 128
    & 16 & 32 & 64 & 128
    & 16 & 32 & 64 & 128
    \\
    \midrule
        
    $\Relax$
    & \multicolumn {4}{c|}{\textbf {\AvgOuter }}
    & \multicolumn {4}{c|}{\textbf {\AvgCPU }}    
    & \multicolumn {4}{c}{\textbf {\AvgInner }}
    \\
    \midrule
     
    \multicolumn {13}{c}{$\Triang^{0} (
    \Ncell=224,
    \coarse{\Ncell}=672)$ \quad \#DOF=(1.4e4, 2.8e4, 5.7e4, 1.1e5)}
    \\
    \midrule
    1  & \cellcolor {gray!20}4 & 4 & 4 & 2 & \cellcolor {gray!20}1.8e-1 & 3.1e-1 & 9.7e-1 & 2.3e1 & \cellcolor {gray!20}2.2 & 2.1 & 2.7 & 1.5\\ 
2e-1 & \cellcolor {gray!20}4 & 4 & 3 & 1 & \cellcolor {gray!20}2.1e-1 & 6.1e-1 & 5.8e0 & 1.9e1 & \cellcolor {gray!20}2.7 & 2.3 & 2.6 & 1.0\\ 
4e-2 & \cellcolor {gray!20}5 & 5 & 4 & 1 & \cellcolor {gray!20}3.3e-1 & 1.1e0 & 7.4e0 & 2.0e1 & \cellcolor {gray!20}1.9 & 2.3 & 1.9 & 1.0\\ 
8e-3 & \cellcolor {gray!20}6 & 4 & 4 & 1 & \cellcolor {gray!20}4.8e-1 & 1.7e0 & 9.0e0 & 1.8e1 & \cellcolor {gray!20}1.6 & 2.0 & 1.8 & 1.0\\ 
2e-3 & \cellcolor {gray!20}5 & 5 & 4 & 2 & \cellcolor {gray!20}5.9e-1 & 1.9e0 & 8.8e0 & 2.0e1 & \cellcolor {gray!20}1.6 & 1.8 & 1.8 & 1.1\\ 
3e-4 & \cellcolor {gray!20}5 & 4 & 4 & 2 & \cellcolor {gray!20}6.4e-1 & 2.1e0 & 8.4e0 & 1.9e1 & \cellcolor {gray!20}1.7 & 2.0 & 1.7 & 1.2\\ 
6e-5 & \cellcolor {gray!20}5 & 4 & 4 & 2 & \cellcolor {gray!20}6.8e-1 & 2.1e0 & 8.8e0 & 2.2e1 & \cellcolor {gray!20}1.7 & 1.8 & 2.0 & 1.2\\ 
1e-5 & \cellcolor {gray!20}4 & 4 & 4 & 2 & \cellcolor {gray!20}7.0e-1 & 2.2e0 & 9.1e0 & 2.4e1 & \cellcolor {gray!20}1.7 & 2.1 & 1.9 & 1.2\\ 
3e-6 & \cellcolor {gray!20}4 & 4 & 4 & 2 & \cellcolor {gray!20}7.1e-1 & 2.4e0 & 7.9e0 & 2.5e1 & \cellcolor {gray!20}1.8 & 1.8 & 1.7 & 1.0\\ 
5e-7 & \cellcolor {gray!20}4 & 4 & 4 & 2 & \cellcolor {gray!20}7.8e-1 & 2.5e0 & 9.2e0 & 2.3e1 & \cellcolor {gray!20}1.7 & 1.9 & 1.9 & 1.1\\ 
\midrule  & \cellcolor {gray!20}4.5 & 4.2 & 4.0 & 1.7 & \cellcolor {gray!20}4.9e-1 & 1.6e0 & 7.4e0 & 2.1e1 & \cellcolor {gray!20}1.9 & 2.0 & 2.0 & 1.1\\ 
\bottomrule 
 \end{tabular} 
      

\begin {tabular}{R{0.8cm }|R{0.45cm }R{0.45cm }R{0.45cm }R{0.52cm}|R{0.8cm}R{0.82cm}R{0.82cm}R{0.7cm}|R{0.7cm }R{0.7cm }R{0.7cm }R{0.7cm }}
\toprule
 
    \multicolumn {13}{c}{$\Triang^{1} (
    \Ncell=896,
    \coarse{\Ncell}=2688)$ \quad \#DOF=(5.6e4, 1.1e5, 2.3e5, 4.6e5)}
    \\
    \midrule
    1  & 4 & \cellcolor {gray!20}4 & 4 & 4 & 5.7e-1 & \cellcolor {gray!20}1.4e0 & 2.6e0 & 1.4e1 & 2.7 & \cellcolor {gray!20}2.3 & 2.1 & 3.5\\ 
2e-1 & 3 & \cellcolor {gray!20}3 & 3 & 3 & 6.5e-1 & \cellcolor {gray!20}1.7e0 & 1.2e1 & 5.8e2 & 2.7 & \cellcolor {gray!20}2.7 & 2.9 & 2.9\\ 
4e-2 & 5 & \cellcolor {gray!20}5 & 5 & 4 & 1.1e0 & \cellcolor {gray!20}3.8e0 & 5.1e1 & 9.0e2 & 2.6 & \cellcolor {gray!20}2.1 & 2.3 & 2.2\\ 
8e-3 & 5 & \cellcolor {gray!20}6 & 5 & 4 & 1.9e0 & \cellcolor {gray!20}8.6e0 & 1.0e2 & 1.1e3 & 3.3 & \cellcolor {gray!20}2.1 & 1.9 & 1.7\\ 
2e-3 & 7 & \cellcolor {gray!20}6 & 5 & 4 & 3.4e0 & \cellcolor {gray!20}1.3e1 & 1.4e2 & 1.0e3 & 6.0 & \cellcolor {gray!20}2.0 & 1.8 & 1.6\\ 
3e-4 & 6 & \cellcolor {gray!20}6 & 5 & 4 & 6.3e0 & \cellcolor {gray!20}1.5e1 & 1.5e2 & 1.0e3 & 14.6 & \cellcolor {gray!20}2.1 & 2.0 & 1.9\\ 
6e-5 & 6 & \cellcolor {gray!20}5 & 4 & 4 & 1.1e1 & \cellcolor {gray!20}1.8e1 & 1.6e2 & 1.0e3 & 32.2 & \cellcolor {gray!20}2.5 & 2.1 & 1.9\\ 
1e-5 & 6 & \cellcolor {gray!20}5 & 4 & 4 & 3.0e1 & \cellcolor {gray!20}1.9e1 & 1.7e2 & 1.1e3 & 95.0 & \cellcolor {gray!20}2.5 & 2.1 & 1.8\\ 
3e-6 & 6 & \cellcolor {gray!20}5 & 4 & 4 & 4.2e1 & \cellcolor {gray!20}2.2e1 & 1.7e2 & 1.0e3 & 122.8 & \cellcolor {gray!20}2.3 & 1.9 & 1.9\\ 
5e-7 & 9 & \cellcolor {gray!20}5 & 4 & 4 & 7.9e1 & \cellcolor {gray!20}2.7e1 & 1.9e2 & 1.1e3 & 151.4 & \cellcolor {gray!20}2.5 & 2.1 & 2.7\\ 
\midrule  & 5.5 & \cellcolor {gray!20}4.8 & 4.3 & 3.7 & 1.5e1 & \cellcolor {gray!20}1.2e1 & 1.0e2 & 8.6e2 & 49.0 & \cellcolor {gray!20}2.3 & 2.1 & 2.2\\ 
\bottomrule 
 \end{tabular} 
      

\begin {tabular}{R{0.8cm }|R{0.45cm }R{0.45cm }R{0.45cm }R{0.52cm}|R{0.8cm}R{0.82cm}R{0.82cm}R{0.7cm}|R{0.7cm }R{0.7cm }R{0.7cm }R{0.7cm }}
\toprule
 
    \multicolumn {13}{c}{$\Triang^{2} (
    \Ncell=3584,
    \coarse{\Ncell}=10752)$ \quad \#DOF=(2.3e5, 4.6e5, 9.1e5, 1.8e6)}
    \\
    \midrule
    1  & 3 & 3 & \cellcolor {gray!20}4 & 4 & 2.3e0 & 5.0e0 & \cellcolor {gray!20}1.0e1 & 2.2e1 & 3.2 & 2.9 & \cellcolor {gray!20}2.4 & 2.3\\ 
2e-1 & 3 & 3 & \cellcolor {gray!20}3 & 3 & 2.9e0 & 6.3e0 & \cellcolor {gray!20}1.5e1 & 2.1e3 & 3.4 & 3.2 & \cellcolor {gray!20}3.1 & 3.3\\ 
4e-2 & 4 & 4 & \cellcolor {gray!20}5 & 9 & 4.9e0 & 2.6e1 & \cellcolor {gray!20}2.1e2 & 1.3e3 & 4.0 & 4.2 & \cellcolor {gray!20}2.5 & 128.5\\ 
8e-3 & 5 & 5 & \cellcolor {gray!20}6 & 20 & 1.1e1 & 9.6e1 & \cellcolor {gray!20}1.2e3 & 2.0e3 & 4.2 & 5.0 & \cellcolor {gray!20}2.4 & 160.1\\ 
2e-3 & 6 & 6 & \cellcolor {gray!20}6 & 55 & 1.9e1 & 1.3e2 & \cellcolor {gray!20}1.5e3 & 6.2e3 & 3.1 & 7.7 & \cellcolor {gray!20}2.0 & 198.8\\ 
3e-4 & 6 & 8 & \cellcolor {gray!20}6 & \dag & 2.5e1 & 2.0e2 & \cellcolor {gray!20}2.3e3 & \dag & 3.3 & 18.1 & \cellcolor {gray!20}2.2 & \dag\\ 
6e-5 & 6 & 6 & \cellcolor {gray!20}4 & \dag & 3.1e1 & 3.8e2 & \cellcolor {gray!20}2.2e3 & \dag & 6.4 & 65.2 & \cellcolor {gray!20}2.6 & \dag\\ 
1e-5 & 6 & 7 & \cellcolor {gray!20}4 & \dag & 4.5e1 & 7.1e2 & \cellcolor {gray!20}1.7e3 & \dag & 13.1 & 133.8 & \cellcolor {gray!20}2.7 & \dag\\ 
3e-6 & 6 & 9 & \cellcolor {gray!20}4 & \dag & 7.3e1 & 1.1e3 & \cellcolor {gray!20}1.7e3 & \dag & 28.9 & 161.5 & \cellcolor {gray!20}2.9 & \dag\\ 
5e-7 & 5 & 16 & \cellcolor {gray!20}4 & \dag & 1.3e2 & 2.1e3 & \cellcolor {gray!20}1.8e3 & \dag & 70.9 & 179.3 & \cellcolor {gray!20}2.6 & \dag\\ 
\midrule  & 4.7 & 6.3 & \cellcolor {gray!20}4.6 & \dag & 2.7e1 & 4.1e2 & \cellcolor {gray!20}1.1e3 & \dag & 12.2 & 79.4 & \cellcolor {gray!20}2.5 & \dag\\ 
\bottomrule 
 \end{tabular} 
      

\begin {tabular}{R{0.8cm }|R{0.45cm }R{0.45cm }R{0.45cm }R{0.52cm}|R{0.8cm}R{0.82cm}R{0.82cm}R{0.7cm}|R{0.7cm }R{0.7cm }R{0.7cm }R{0.7cm }}
\toprule
 
    \multicolumn {13}{c}{$\Triang^{3} (
    \Ncell=14336,
    \coarse{\Ncell}=43008)$ \quad \#DOF=(9.0e5, 1.8e6, 3.7e6, 7.3e6)}
    \\
    \midrule
    1  & 3 & 3 & 3 & \cellcolor {gray!20}3 & 9.3e0 & 2.0e1 & 5.0e1 & \cellcolor {gray!20}9.8e1 & 3.5 & 3.4 & 2.9 & \cellcolor {gray!20}2.8\\ 
2e-1 & 2 & 3 & 3 & \cellcolor {gray!20}3 & 1.1e1 & 2.5e1 & 6.0e1 & \cellcolor {gray!20}4.7e2 & 4.0 & 3.9 & 3.6 & \cellcolor {gray!20}3.3\\ 
4e-2 & 4 & 4 & 4 & \cellcolor {gray!20}7 & 2.7e1 & 1.5e2 & 3.4e3 & \cellcolor {gray!20}1.8e3 & 5.6 & 4.6 & 22.9 & \cellcolor {gray!20}84.2\\ 
8e-3 & 5 & 5 & 9 & \cellcolor {gray!20}12 & 7.4e1 & 5.0e2 & 2.6e3 & \cellcolor {gray!20}3.7e3 & 5.7 & 5.2 & 154.9 & \cellcolor {gray!20}105.4\\ 
2e-3 & 6 & 7 & 18 & \cellcolor {gray!20}26 & 1.4e2 & 2.4e3 & 4.9e3 & \cellcolor {gray!20}1.1e4 & 4.5 & 3.3 & 182.5 & \cellcolor {gray!20}152.1\\ 
3e-4 & 6 & 6 & 58 & \cellcolor {gray!20}59 & 2.1e2 & 2.7e3 & 1.5e4 & \cellcolor {gray!20}3.6e4 & 4.0 & 3.2 & 199.0 & \cellcolor {gray!20}198.9\\ 
6e-5 & 5 & 6 & \dag & \cellcolor {gray!20}162 & 2.6e2 & 2.1e3 & \dag & \cellcolor {gray!20}9.5e4 & 4.3 & 4.1 & \dag & \cellcolor {gray!20}200.0\\ 
1e-5 & 5 & 5 & \dag & \cellcolor {gray!20}\dag & 3.2e2 & 2.4e3 & \dag & \cellcolor {gray!20}\dag & 4.5 & 7.1 & \dag & \cellcolor {gray!20}\dag\\ 
3e-6 & 5 & 5 & \dag & \cellcolor {gray!20}\dag & 3.5e2 & 2.7e3 & \dag & \cellcolor {gray!20}\dag & 4.4 & 14.9 & \dag & \cellcolor {gray!20}\dag\\ 
5e-7 & 5 & 5 & \dag & \cellcolor {gray!20}\dag & 3.6e2 & 3.7e3 & \dag & \cellcolor {gray!20}\dag & 4.6 & 41.8 & \dag & \cellcolor {gray!20}\dag\\ 
\midrule  & 4.3 & 4.7 & \dag & \cellcolor {gray!20}\dag & 1.5e2 & 1.4e3 & \dag & \cellcolor {gray!20}\dag & 4.6 & 7.9 & \dag & \cellcolor {gray!20}\dag\\ 
\bottomrule 
 \end{tabular} 
    }
  \end{center}
   \caption{Numerical results using the preconditioner based on the primal Schur complement. Each sub-table reports the results for a given mesh, from the coarsest (top) to the finest (bottom), and for different time discretizations ($\Deltat=1/(\Ntime+1)$). Each sub-table reports, while $\Relax$ is reduced (leftmost column), the averaged outer iteration and CPU time per linear system (\AvgOuter{} and \AvgCPU{}) and the averaged inner iterations per outer iteration \AvgInner{} for solving the linear system in~\cref{eq:primal-linsys} (metrics defined in~\cref{sec:preconditioning}). A final row summarizes the averages on the whole simulation. We highlighted in gray the time-space combination providing the uniform discretization. The $\dag$ symbol denotes those \IP{} steps where the linear solver failed.
   } 
  \label{tab:sol10perc}
\end{table}
\addtolength{\tabcolsep}{3.5pt}
In~\cref{tab:sol10perc} we summarize the results obtained using the preconditioner presented in this section for different combinations of time and space discretizations. We use the metrics described in~\cref{sec:preconditioning}, where the average number of inner iterations, denoted by \AvgInner{}, refers to the linear system in~\cref{eq:primal-linsys}.

The number of average outer iterations, denoted by \AvgOuter{}, remains practically constant only for the coarsest and most uniform combinations of grids and time steps. However, in some cases, also \AvgOuter{} increases when the average number of inner iterations \AvgInner{} get closer to 200 (the limit we fixed for the inner solver), which means that, in some cases, the inner solver did not reach the accuracy $\InnerTol=1e-1$. For the cases with more degrees of freedom, this led to a huge CPU time (almost $1e5$ seconds per linear system) and to failures occurring as $\Relax$ is reduced.

The average number of inner iterations \AvgInner{} is strongly influenced by the relationship between the mesh and the time step $\Deltat=1/(\Ntimestep+1)$.  The phenomenon is particularly evident for mesh $\Triang^2$. For $\Ntimestep+1=32$, the average number of inner iterations increases progressively while $\Relax$ is reduced, reaching $\approx 179$ at the last IP iteration. For $\Ntimestep+1=128$, \AvgInner{} increases reaching the 200 iterations limit. For $\Ntimestep+1=64$ instead, which corresponds to the most balanced time-space scaling, this number remains between $2.0$ and $2.7$. The same phenomenon occurs using $\Triang^1$, passing from $\Ntimestep+1=16$ to $\Ntimestep+1=32$. We attribute this to an improper construction of the coarse matrix sequence used by the multigrid solver due to the poor scaling of the temporal and spatial components of $\Matr{\Sblock}_p$. The AGMG is based on an aggregation-based coarsening approach, hence if these components are unbalanced the coarse matrices may not reflect the infinite-dimensional operators discretized by our problem, compromising the performance of the multigrid solver.

This phenomenon can only get worse in the last IP iterations due to presence of the term $\Cblock^{-1}=\Matr{\BlockMassTdens}^{-1}\Diag{\Vect{\Tdens}}\Diag{\Vect{\Slack}}^{-1}$ in $\Sblock_p=\Ablock+\Bblock \Cblock^{-1}\Bblock^T$. In fact, the term $\Cblock^{-1}$ contains entries varying by several order of magnitude due to the relaxed complementarity condition $\Vect{\Slack}\odot \Vect{\Tdens}\approx\mu \Vect{1}$. This also means that, on those entries where $\Vect{\Tdens}$ remain strictly positive, the anisotropic Laplacian $\Matr{\Sblock}_{xx}$ scales like $1/\Relax$, making the linear system~\cref{eq:primal-linsys} more difficult to solve since this anisotropic term becomes the dominant part of $\Matr{\Sblock}_p$. Similar matrices and analogous considerations can be found in~\cite{behata:2009,benzi2011preconditioning,haber2015multilevel,morini2017comparison}.

Furthermore, even in best-case scenarios, we experimented that the preprocessing phase is the most demanding in terms of CPU time. Typically, more than 50\% of the time is spent building the sequence of coarse matrices, reaching peaks of $90\%$ during the last \IP{} steps using the finest discretization.

In order to cope with all the issues presented so far, we tried to neglect some components of $\Matr{\Sblock}_p$. We tried to remove the extra-diagonal or lower-diagonal part of $\Sblock_p$, or the time-space operator $\Matr{\Sblock}_{tx}+\Matr{\Sblock}_{tx}^T$ so that the remaining terms form a weighted time-space Laplacian (similarly to what done in~\cite{haber2015multilevel}). However, none of these approaches worked. All numerical experiments suggest that none of these terms can be neglected without altering the spectral property of the matrix, and consequently affecting the effectiveness of the preconditioner.

Summarizing the results obtained in this section, the preconditioner based on the primal Schur complement $\Sblock_p=\Ablock + \Bblock^T \Cblock^{-1} \Bblock$ lacks both in efficiency and in robustness, due to the complexity of the PDE behind the linear system in~\cref{eq:primal-linsys}, in particular when the relaxation parameter $\Relax$ goes to zero.

\subsection{The SIMPLE preconditioner}
\label{sec:dual-schur}
In this section, we recall the classical SIMPLE preconditioner described in~\cite{Patankar72,Patankar80} and based on the approximation of the inverse of the dual Schur complement $\Matr{\Sblock}=-(\Cblock+\Bblock{\Ablock}^{-1} \Bblock^T)$ (we recall that, formally, we cannot write ${\Ablock}^{-1}$ since the matrix $\Ablock$ is singular). The SIMPLE preconditioner is given by
\begin{gather*}
  \Matr{P}^{-1}
  =
  \begin{pmatrix}
    \Matr{\Identity{}}& -\Matr{\Approx{\Amat}}^{-1}\Bblock^T 
    \\
    \ZeroMat& \Matr{\Identity{}}
  \end{pmatrix}
  \begin{pmatrix}
    \Matr{\Approx{\Amat}}^{-1} & \ZeroMat
    \\
    \ZeroMat & \Matr{\Approx{\Smat}}^{-1}
  \end{pmatrix}
  \begin{pmatrix}
    \Matr{\Identity{}} & \ZeroMat
    \\
    -\Bblock\Matr{\Approx{\Amat}}^{-1} & \Matr{\Identity{}}
  \end{pmatrix}
\end{gather*}
where
\begin{equation*}
  \Matr{\Approx{\Amat}}=\Diag{\Ablock},
  \quad
  \Matr{\Approx{\Smat}}=-(\Cblock+\Bblock\Matr{\Approx{\Amat}}^{-1}\Bblock^T)\ .
\end{equation*}
This preconditioner is referred to as constraint preconditioner in the literature studying linear solvers for \IP{} methods~\cite{morini2017comparison}.
The main computational cost of applying this preconditioner is to solve linear systems in the form
\begin{equation}
  \label{eq:linsys-simple}
  \Matr{\Approx{\Smat}} \Vect{y}= \Vect{c} \, .
\end{equation}
The matrix $\Matr{\Approx{\Smat}}$ can be formed explicitly and is block tridiagonal.

During the latest steps of the \IP{} method, when $\Relax\to0$, linear system~\cref{eq:linsys-simple} may become ill-conditioned since the diagonal matrix $\Cblock=\Matr{\BlockMassTdens}\Diag{\Vect{\Tdens}}^{-1}\Diag{\Vect{\Slack}}$ can contain terms that vary by several orders of magnitude, since $\Vect{\Tdens}\odot \Vect{\Slack}\approx \Relax \Vect{1}\to \Vect{0}$.  To cope with this issue, we scale both sides of~\cref{eq:linsys-simple} by $\Diag{\Vect{\Tdens}}$ and we solve the linear system
\begin{equation} \label{eq:linsys-simple-mode2}
  -(\MassTdensBlock\Diag{\Vect{\Slack}} + \Diag{\Vect{\Tdens}}\Bblock\Matr{\Approx{\Amat}}^{-1}\Bblock^T) \Vect{y}=     \Diag{\Vect{\Tdens}}\Vect{c}\,.
\end{equation}
This scaling may seem a dangerous procedure if small entries appears in $\Vect{\Tdens}$, however it is used within the preconditioner application, where the potential inaccuracies do not affect the convergence of FGMRES to the solution of the linear system.
We use the AGMG solver with tolerance $\InnerTol=1e-1$ to solve the linear system in~\cref{eq:linsys-simple-mode2} and denote by $\Matr{\Approx{\Smat}}^{-1}_{\InnerTol}$ the resulting (nonlinear) operator. We found experimentally that this value provided the best performance. Further reductions of $\InnerTol$ led only to an increase of the number of inner iterations but no reduction in the outer loop iterations.

\addtolength{\tabcolsep}{-3.5pt}
\begin{table}
  \begin{center}
    {\footnotesize
      

\begin {tabular}{R{0.65cm }|R{0.6cm }R{0.6cm }R{0.6cm }R{0.6cm}|R{0.82cm }R{0.82cm }R{0.82cm }R{0.82cm }|R{0.4cm }R{0.55cm }R{0.55cm }R{0.55cm }}
\toprule
$\Ntime\!+\!1$&16 & 32 & 64 & 128&16 & 32 & 64 & 128&16 & 32 & 64 & 128\\\midrule
    $\Relax$
    & \multicolumn {4}{c|}{\textbf {\AvgOuter}}
    & \multicolumn {4}{c|}{\textbf {\AvgCPU }}
    & \multicolumn {4}{c}{\textbf {\AvgInner }}
    \\
    \midrule
     
    \multicolumn {13}{c}{$\Triang^{0} (
    \Ncell=224,
    \coarse{\Ncell}=672)$ \quad \#DOF=(1.4e4, 2.8e4, 5.7e4, 1.1e5)}
    \\
    \midrule
    1  & \cellcolor {gray!20}108 & 110 & 106 & 107 & \cellcolor {gray!20}3.8e-1 & 7.3e-1 & 1.6e0 & 4.1e0 & \cellcolor {gray!20}1.0 & 1.0 & 1.8 & 2.1\\ 
2e-1 & \cellcolor {gray!20}72 & 72 & 73 & 74 & \cellcolor {gray!20}3.1e-1 & 6.3e-1 & 1.4e0 & 3.4e0 & \cellcolor {gray!20}1.1 & 1.6 & 2.4 & 2.9\\ 
4e-2 & \cellcolor {gray!20}65 & 65 & 66 & 66 & \cellcolor {gray!20}3.1e-1 & 6.4e-1 & 1.5e0 & 3.5e0 & \cellcolor {gray!20}1.7 & 2.7 & 3.5 & 4.2\\ 
8e-3 & \cellcolor {gray!20}56 & 58 & 60 & 60 & \cellcolor {gray!20}2.9e-1 & 6.4e-1 & 1.5e0 & 3.8e0 & \cellcolor {gray!20}2.2 & 3.1 & 4.0 & 4.8\\ 
2e-3 & \cellcolor {gray!20}52 & 54 & 48 & 52 & \cellcolor {gray!20}2.9e-1 & 6.4e-1 & 1.4e0 & 3.7e0 & \cellcolor {gray!20}2.3 & 3.1 & 4.2 & 5.2\\ 
3e-4 & \cellcolor {gray!20}52 & 43 & 47 & 48 & \cellcolor {gray!20}3.1e-1 & 5.8e-1 & 1.4e0 & 3.8e0 & \cellcolor {gray!20}2.5 & 3.4 & 4.4 & 6.5\\ 
6e-5 & \cellcolor {gray!20}42 & 41 & 41 & 39 & \cellcolor {gray!20}3.0e-1 & 5.8e-1 & 1.3e0 & 3.1e0 & \cellcolor {gray!20}2.9 & 4.4 & 5.6 & 7.3\\ 
1e-5 & \cellcolor {gray!20}41 & 42 & 40 & 43 & \cellcolor {gray!20}2.8e-1 & 6.3e-1 & 1.3e0 & 3.8e0 & \cellcolor {gray!20}3.8 & 5.5 & 6.6 & 8.5\\ 
3e-6 & \cellcolor {gray!20}43 & 55 & 50 & 68 & \cellcolor {gray!20}3.2e-1 & 2.0e0 & 1.8e0 & 6.5e0 & \cellcolor {gray!20}3.8 & 19.5 & 7.9 & 9.7\\ 
5e-7 & \cellcolor {gray!20}54 & 72 & 67 & 74 & \cellcolor {gray!20}3.8e-1 & 2.0e0 & 2.7e0 & 9.7e0 & \cellcolor {gray!20}3.8 & 14.0 & 8.7 & 14.1\\ 
\toprule  & \cellcolor {gray!20}60 & 63 & 61 & 65 & \cellcolor {gray!20}3.2e-1 & 8.6e-1 & 1.6e0 & 4.4e0 & \cellcolor {gray!20}2.1 & 5.0 & 4.3 & 5.7\\ 
\bottomrule 
 \end{tabular}
      

\begin {tabular}{R{0.65cm }|R{0.6cm }R{0.6cm }R{0.6cm }R{0.6cm}|R{0.82cm }R{0.82cm }R{0.82cm }R{0.82cm }|R{0.4cm }R{0.55cm }R{0.55cm }R{0.55cm }}
\toprule
 
    \multicolumn {13}{c}{$\Triang^{1} (
    \Ncell=896,
    \coarse{\Ncell}=2688)$ \quad \#DOF=(5.6e4, 1.1e5, 2.3e5, 4.6e5)}
    \\
    \midrule
    1  & 250 & \cellcolor {gray!20}230 & 229 & 230 & 3.1e0 & \cellcolor {gray!20}5.9e0 & 1.2e1 & 2.8e1 & 1.0 & \cellcolor {gray!20}1.0 & 1.1 & 1.9\\ 
2e-1 & 136 & \cellcolor {gray!20}144 & 147 & 148 & 1.9e0 & \cellcolor {gray!20}4.3e0 & 9.8e0 & 2.3e1 & 1.0 & \cellcolor {gray!20}1.2 & 1.6 & 2.4\\ 
4e-2 & 114 & \cellcolor {gray!20}116 & 119 & 119 & 1.7e0 & \cellcolor {gray!20}4.0e0 & 9.7e0 & 2.3e1 & 1.4 & \cellcolor {gray!20}1.8 & 3.0 & 3.9\\ 
8e-3 & 90 & \cellcolor {gray!20}107 & 107 & 106 & 1.5e0 & \cellcolor {gray!20}4.0e0 & 9.9e0 & 2.4e1 & 1.8 & \cellcolor {gray!20}2.2 & 3.8 & 4.7\\ 
2e-3 & 89 & \cellcolor {gray!20}91 & 89 & 91 & 1.6e0 & \cellcolor {gray!20}3.8e0 & 9.5e0 & 2.4e1 & 2.0 & \cellcolor {gray!20}2.9 & 4.3 & 5.3\\ 
3e-4 & 102 & \cellcolor {gray!20}85 & 84 & 79 & 1.9e0 & \cellcolor {gray!20}3.8e0 & 9.7e0 & 2.3e1 & 2.0 & \cellcolor {gray!20}3.0 & 4.5 & 5.4\\ 
6e-5 & 105 & \cellcolor {gray!20}97 & 82 & 74 & 2.0e0 & \cellcolor {gray!20}4.7e0 & 1.0e1 & 2.3e1 & 2.0 & \cellcolor {gray!20}3.1 & 4.8 & 5.9\\ 
1e-5 & 100 & \cellcolor {gray!20}86 & 79 & 78 & 2.1e0 & \cellcolor {gray!20}4.2e0 & 1.1e1 & 2.6e1 & 2.2 & \cellcolor {gray!20}3.2 & 6.1 & 7.1\\ 
3e-6 & 91 & \cellcolor {gray!20}74 & 78 & 85 & 2.1e0 & \cellcolor {gray!20}4.2e0 & 1.2e1 & 3.4e1 & 2.7 & \cellcolor {gray!20}4.7 & 7.9 & 10.2\\ 
5e-7 & 88 & \cellcolor {gray!20}74 & 96 & 115 & 2.1e0 & \cellcolor {gray!20}4.7e0 & 1.6e1 & 5.2e1 & 3.1 & \cellcolor {gray!20}5.3 & 8.3 & 12.3\\ 
\toprule  & 120 & \cellcolor {gray!20}118 & 118 & 119 & 2.0e0 & \cellcolor {gray!20}4.4e0 & 1.1e1 & 2.7e1 & 1.7 & \cellcolor {gray!20}2.2 & 3.4 & 4.7\\ 
\bottomrule 
 \end{tabular}
      

\begin {tabular}{R{0.65cm }|R{0.6cm }R{0.6cm }R{0.6cm }R{0.6cm}|R{0.82cm }R{0.82cm }R{0.82cm }R{0.82cm }|R{0.4cm }R{0.55cm }R{0.55cm }R{0.55cm }}
\toprule
 
    \multicolumn {13}{c}{$\Triang^{2} (
    \Ncell=3584,
    \coarse{\Ncell}=10752)$ \quad \#DOF=(2.3e5, 4.6e5, 9.1e5, 1.8e6)}
    \\
    \midrule
    1  & 516 & 459 & \cellcolor {gray!20}497 & 475 & 2.6e1 & 4.8e1 & \cellcolor {gray!20}1.2e2 & 2.7e2 & 1.0 & 1.0 & \cellcolor {gray!20}1.0 & 1.0\\ 
2e-1 & 300 & 298 & \cellcolor {gray!20}319 & 317 & 1.7e1 & 3.5e1 & \cellcolor {gray!20}8.4e1 & 2.1e2 & 1.0 & 1.0 & \cellcolor {gray!20}1.1 & 1.6\\ 
4e-2 & 206 & 222 & \cellcolor {gray!20}234 & 221 & 1.3e1 & 2.9e1 & \cellcolor {gray!20}7.2e1 & 1.8e2 & 1.5 & 1.4 & \cellcolor {gray!20}1.9 & 3.0\\ 
8e-3 & 190 & 164 & \cellcolor {gray!20}167 & 169 & 1.3e1 & 2.4e1 & \cellcolor {gray!20}5.8e1 & 1.5e2 & 1.8 & 1.8 & \cellcolor {gray!20}2.4 & 4.0\\ 
2e-3 & 185 & 161 & \cellcolor {gray!20}150 & 154 & 1.4e1 & 2.5e1 & \cellcolor {gray!20}5.6e1 & 1.6e2 & 2.6 & 2.0 & \cellcolor {gray!20}3.2 & 4.8\\ 
3e-4 & 223 & 170 & \cellcolor {gray!20}168 & 139 & 1.8e1 & 2.8e1 & \cellcolor {gray!20}6.9e1 & 1.6e2 & 3.0 & 2.5 & \cellcolor {gray!20}3.5 & 5.2\\ 
6e-5 & 228 & 222 & \cellcolor {gray!20}157 & 141 & 2.0e1 & 3.8e1 & \cellcolor {gray!20}7.1e1 & 1.8e2 & 3.0 & 2.6 & \cellcolor {gray!20}3.7 & 5.5\\ 
1e-5 & 255 & 183 & \cellcolor {gray!20}156 & 145 & 2.3e1 & 3.4e1 & \cellcolor {gray!20}7.4e1 & 1.9e2 & 3.0 & 2.7 & \cellcolor {gray!20}3.9 & 5.7\\ 
3e-6 & 250 & 189 & \cellcolor {gray!20}162 & 136 & 2.4e1 & 3.8e1 & \cellcolor {gray!20}8.1e1 & 1.9e2 & 3.1 & 2.9 & \cellcolor {gray!20}4.1 & 6.2\\ 
5e-7 & 252 & 160 & \cellcolor {gray!20}140 & 140 & 2.5e1 & 3.3e1 & \cellcolor {gray!20}7.4e1 & 2.2e2 & 3.1 & 3.0 & \cellcolor {gray!20}4.2 & 6.9\\ 
\toprule  & 272 & 236 & \cellcolor {gray!20}237 & 224 & 1.9e1 & 3.4e1 & \cellcolor {gray!20}7.8e1 & 2.0e2 & 1.9 & 1.7 & \cellcolor {gray!20}2.1 & 3.1\\ 
\bottomrule 
 \end{tabular}
      

\begin {tabular}{R{0.65cm }|R{0.6cm }R{0.6cm }R{0.6cm }R{0.6cm}|R{0.82cm }R{0.82cm }R{0.82cm }R{0.82cm }|R{0.4cm }R{0.55cm }R{0.55cm }R{0.55cm }}
\toprule
 
    \multicolumn {13}{c}{$\Triang^{3} (
    \Ncell=14336,
    \coarse{\Ncell}=43008)$ \quad \#DOF=(9.0e5, 1.8e6, 3.7e6, 7.3e6)}
    \\
    \midrule
    1  & 1003 & 919 & 959 & \cellcolor {gray!20}1045 & 2.2e2 & 4.9e2 & 1.1e3 & \cellcolor {gray!20}3.8e3 & 1.0 & 1.0 & 1.0 & \cellcolor {gray!20}1.0\\ 
2e-1 & 1124 & \dag & \dag & \cellcolor {gray!20}\dag & 2.7e2 & \dag & \dag & \cellcolor {gray!20}\dag & 1.0 & \dag & \dag & \cellcolor {gray!20}\dag\\ 
4e-2 & \dag & \dag & 540 & \cellcolor {gray!20}503 & \dag & \dag & 6.9e2 & \cellcolor {gray!20}2.0e3 & \dag & \dag & 1.4 & \cellcolor {gray!20}1.8\\ 
8e-3 & 368 & 388 & 404 & \cellcolor {gray!20}347 & 1.1e2 & 2.5e2 & 5.6e2 & \cellcolor {gray!20}1.5e3 & 1.8 & 1.8 & 1.8 & \cellcolor {gray!20}2.4\\ 
2e-3 & 397 & 313 & 279 & \cellcolor {gray!20}318 & 1.6e2 & 2.5e2 & 4.2e2 & \cellcolor {gray!20}1.5e3 & 2.8 & 2.8 & 2.1 & \cellcolor {gray!20}3.5\\ 
3e-4 & 578 & 348 & 297 & \cellcolor {gray!20}290 & 3.2e2 & 2.9e2 & 4.8e2 & \cellcolor {gray!20}1.5e3 & 3.3 & 3.0 & 2.7 & \cellcolor {gray!20}3.9\\ 
6e-5 & 419 & 483 & 327 & \cellcolor {gray!20}319 & 3.1e2 & 4.3e2 & 5.7e2 & \cellcolor {gray!20}1.7e3 & 3.9 & 3.3 & 3.0 & \cellcolor {gray!20}4.0\\ 
1e-5 & 407 & 364 & 395 & \cellcolor {gray!20}284 & 3.6e2 & 3.6e2 & 7.2e2 & \cellcolor {gray!20}1.6e3 & 4.2 & 3.8 & 3.0 & \cellcolor {gray!20}4.2\\ 
3e-6 & 389 & 301 & 318 & \cellcolor {gray!20}251 & 3.8e2 & 3.1e2 & 6.0e2 & \cellcolor {gray!20}1.5e3 & 4.4 & 3.9 & 3.0 & \cellcolor {gray!20}4.6\\ 
5e-7 & 410 & 305 & 411 & \cellcolor {gray!20}253 & 3.5e2 & 3.3e2 & 8.1e2 & \cellcolor {gray!20}1.6e3 & 4.9 & 4.0 & 3.0 & \cellcolor {gray!20}5.0\\ 
\toprule  & \dag & \dag & \dag & \cellcolor {gray!20}\dag & \dag & \dag & \dag & \cellcolor {gray!20}\dag & \dag & \dag & \dag & \cellcolor {gray!20}\dag\\ 
\bottomrule 
 \end{tabular}
    }
  \end{center}
  \caption{
    Numerical results using the SIMPLE preconditioner. Each sub-table reports the results for a given mesh, from the coarsest (top) to the finest (bottom), and for different time discretizations ($\Deltat=1/(\Ntime+1)$). Each sub-table reports, while $\Relax$ is reduced (leftmost column), the averaged outer iteration and CPU time per linear system (\AvgOuter{} and \AvgCPU{}) and the averaged inner iterations per outer iteration \AvgInner{} for solving the linear system in~\cref{eq:linsys-simple-mode2}  (metrics defined in~\cref{sec:preconditioning}). A final row summarizes the averages on the whole simulation. We highlighted in gray the time-space combination providing the uniform discretization. The $\dag$ symbol denotes those \IP{} steps where the linear solver failed.
  }
  \label{tab:sol11}
\end{table}
\addtolength{\tabcolsep}{3.5pt}
The results for the SIMPLE preconditioner are summarized in~\cref{tab:sol11}. Despite its simplicity, the proposed preconditioner turned out to be robust (in particular for small values of $\Relax$) but not particularly efficient. Some failures occurred only at the initial \IP{} steps using the finest grid. In this cases, we restarted the solver using the data obtained with the $\BB$-preconditioner in order to study the behavior of the SIMPLE preconditioner close to the optimal solution. The average number of inner iterations \AvgInner{} increases only sightly as $\Relax \to 0$, independently of the time step and the mesh size used or the \IP{} step considered. The preprocessing time required by the AGMG solver is limited, approximately accounting for $1\%$ of the CPU time spent in the linear system solution. The average number of outer iterations \AvgOuter{} is affected only mildly by the number of time steps $\Ntime+1$. 

Remarkably, the preconditioner becomes more efficient as the \IP{} relaxation term $\mu$ goes to zero. However, \AvgOuter{} more than doubles at each mesh refinement. Both phenomena can also be explained by looking at the structure of the preconditioned matrix, which is given by
\begin{align*}
   \begin{pmatrix}
    \Ablock
    &
    \Bblock^T
    \\
    \Bblock
    &
    - \Cblock
  \end{pmatrix}
  \Matr{P}^{-1}
  &\!=\!
  \begin{pmatrix}
    \Matr{\Identity{}}
    +
    (\Ablock\Matr{\Approx{\Amat}}^{-1}\!\!\!\!-\Matr{\Identity{}})
    (
    \Matr{\Identity{}}
    \!+
    \Bblock^T \Matr{\Approx{\Smat}}^{-1}_{\InnerTol}
    \Bblock
    \Matr{\Approx{\Amat}}^{-1}
    )
    &
    (\Matr{\Identity{}}\!\!\!-\Ablock\Matr{\Approx{\Amat}}^{-1})
    \Bblock^T \Matr{\Approx{\Smat}}^{-1}_{\InnerTol}
    \\
    (\Matr{\Identity{}} - \Matr{\Approx{\Smat}}\Matr{\Approx{\Smat}}^{-1}_{\InnerTol})\Bblock\Matr{\Approx{\Amat}}^{-1}
    & \Matr{\Approx{\Smat}}\Matr{\Approx{\Smat}}^{-1}_{\InnerTol}
  \end{pmatrix},
\end{align*}
and estimating its eigenvalues. In order to do this, let us first observe that the preconditioned matrix is approximately a block upper triangular matrix, whose bottom-right block becomes close to the identity plus a perturbation of order $\InnerTol=1e-1$, hence there are $\Nr$ eigenvalues close to one.

The upper left block is the identity plus a perturbation matrix (which should be as small as possible to have an ideal preconditioner) given by product of two factors. The first factor is 
  \[
  \Ablock\Matr{\Approx{\Amat}}^{-1}-\Matr{\Identity{}}
  =\BlkDiag{
    \Matr{A}^{\tstep}\Diag{\Matr{A}^{\tstep}}^{-1}-\Identity{\fine{\Ncell}}
  }.
  \]
  Since the matrices $(\Matr{A}^{\tstep})_{\tstep=1,\ldots,\Ntime+1}$ are weighted Laplacian matrices, a diagonal preconditioner becomes a poor preconditioner for large $\Matr{A}^{\tstep}$ (see~\cite{kamenski2014conditioning} for precise estimates). Hence the largest eigenvalue of the first factor will increase while refining in space, but not in time, by the block diagonal structure of matrix $\Ablock$.
The second factor, which we would like to keep as small as possible to compensate the first one, is $\Identity{}+\Matr{Q}$  with $\Matr{Q}=\Bblock^T \Matr{\Approx{\Smat}}^{-1}_{\InnerTol} \Bblock \Matr{\Approx{\Amat}}^{-1}$.
Using the formula in~\cref{eq:linsys-simple-mode2} and the complementarity condition $\Vect{\Tdens}\odot\Vect{\Slack}\approx \Relax \Vect{1}$, we can approximate it as follows:
  \begin{align*}
  \Identity{}+\Matr{Q}
  &=
  \Identity{}-\Bblock^T
  (\MassTdensBlock\Diag{\Vect{\Slack}} + \Diag{\Vect{\Tdens}}\Bblock\Matr{\Approx{\Amat}}^{-1}\Bblock^T)^{-1}
  \Diag{\Vect{\Tdens}}  \Bblock \Matr{\Approx{\Amat}}^{-1}
  \\
    &\approx
  \Identity{}-\Bblock^T
  (\MassTdensBlock\Relax + \Diag{\Vect{\Tdens}}^2\Bblock\Matr{\Approx{\Amat}}^{-1}\Bblock^T)^{-1}
  \Diag{\Vect{\Tdens}}^2  \Bblock \Matr{\Approx{\Amat}}^{-1}.
  \end{align*}
  If we set $\mu=0$ in this formula, the resulting matrix has 0 or 1 eigenvalues since it is idempotent. This gives a qualitative explanation of why the SIMPLE preconditioner performs better as $\Relax \to 0$. Nevertheless, this approach is not viable to tackle large problems due to the quadratic complexity with respect to the number of spatial degrees of freedom.

\subsection{The $\BB$-preconditioner}
\label{sec:commute}
We present now the third preconditioner considered in this paper. Its design started from the idea of using the following block triangular preconditioner (see for example~\cite{Benzi-et-al:2005}),
\begin{equation*}
  \Matr{P}_{t}
  =
  \begin{pmatrix}
    \Matr{\Amat}& \Bblock^T
    \\
    \ZeroMat&  -\Matr{\Approx{\Smat}}
  \end{pmatrix},
\end{equation*}
where $\Matr{\Approx{\Smat}}$ should ideally equal the dual Schur complement $\Sblock=-(\Cblock+\Bblock{\Ablock}^{-1} \Bblock^T)$ or some approximation of it (we again recall that, formally, we cannot write ${\Ablock}^{-1}$ since the matrix $\Ablock$ is singular). Given a vector $(\Vect{\RhsPot};\Vect{\RhsTdens})\in \REAL^{\Np+\Nr}$, one application of this preconditioner requires to compute the vector $(\Vect{x};\Vect{y})\in \REAL^{\Np+\Nr}$ solving
\begin{align}
  \label{eq:wls-wrong}
  -\Matr{\Approx{\Smat}} \Vect{y} &= \Vect{\RhsTdens},
  \\
  \label{eq:x-triang-wrong}
  \Ablock \Vect{x} & = \Vect{b} = \Vect{\RhsPot} - \Bblock^T \Vect{y}.
\end{align}
This approach is particularly attractive from the computational point of view, since matrix $\Ablock$ in~\cref{eq:x-triang} is block diagonal and each block $\Block{\tstep}{\Matr{A}}$ is a weighted Laplacian. Thus, we can efficiently get the solution $\Vect{x}$ by solving separately $\Ntimep$ linear systems in the form $\Block{\tstep}{\Matr{A}} \Block{\tstep}{\Vect{x}}=\Block{\tstep}{\Vect{b}}$ using multigrid solvers. However, two issues arise.

The first one is due to the fact that each block $\Block{\tstep}{\Matr{A}}$ is singular, thus we need to ensure that 
\begin{equation}
  \label{eq:ortogonality}
\Block{\tstep}{\Vect{b}}\in \Image{\Block{\tstep}{\Matr{A}}}=\Kernel{\Block{\tstep}{\Matr{A}}}^\perp=\left\{\Vect{1}\in \REAL^{\fine{\Ncell}}\right\}^\perp,
\quad \tstep=1,\ldots\Ntimep\,,
\end{equation}
for the linear systems to be admissible. This condition may not hold for the right-hand side $\Vect{b}=\Vect{c}-\Bblock^T\Vect{y}$ in~\cref{eq:x-triang-wrong} for any application of the preconditioner, because $\Image{\Bblock^T}$ is not orthogonal to $\Kernel{\Ablock}$ (see~\cref{remark:kernel}).
Even if the orthogonality condition in~\cref{eq:ortogonality} is satisfied, or imposed by orthogonalization, there are multiple solutions, which differ by an additive constant, and thus we also need a criteria to select one. This is not trivial since the solution $\Vect{\DPot}$ of~\cref{eq:saddle-point} must also solve the second block of equations, $\Bblock \Vect{\DPot}-\Cblock \Vect{\DTdens} =\Vect{\tilde{g}}$. In the next section, we describe how to deal with this compatibility issue with a reformulation of the continuous nonlinear system in~\cref{eq:l2}.

A second, harder, issue is the design of a good approximation of the inverse of the matrix $\Matr{\Smat}=-(\Cblock+\Bblock \Matr{{\Amat}}^{-1} \Bblock^T)$, used in the solution of the linear system~\cref{eq:wls-wrong}. First, $\Matr{{\Smat}}$ is not well defined since $\Ablock$ is singular and $\Image{\Bblock^T}$ is not orthogonal to $\Kernel{\Ablock}$. Moreover, even overcoming this problem, the matrix $\Matr{\Smat}$ is block tridiagonal but with dense blocks, due to the presence of the pseudo-inverse of the Laplacian matrices $\Block{\tstep}{\Matr{A}}$, thus too costly to form and invert. This second issue is addressed in~\cref{sec:bb_definition}.

\subsubsection{Splitting formulation}
In order to cope with the first issue, we split the solution $\Pot$ of~\cref{eq:l2} into two components, a new potential $\TPot:\Domain\times[0,1]\to \REAL$ satisfying the constraint $\int_{\Domain}\TPot(t)\dx=0$ for all $t\in[0,1]$ and a correction term:
\begin{equation}
  \label{eq:pot-shift}
  \Pot(t,x) = \TPot(t,x) + \int_{0}^{t} \Shift(s)ds
\end{equation}
where $\Shift:[0,1]\rightarrow \REAL$ depends only on the time variable.
In the new unknowns $(\TPot,\Tdens,\Slack,\Shift)$ the system of PDEs \eqref{eq:l2} becomes
\begin{equation}
  \label{eq:l2-shift}
  \begin{aligned}
    &-\Dt{\Tdens}-\Div_{x}\left(\Tdens \Grad_{x} \TPot \right) = 0, \\
    &\Dt{\TPot} + \frac{|\Grad_{x} \TPot|^2}{2}  + \Slack + \Shift = 0, \\  
    &\Tdens \geq 0, \, \Slack\geq 0, \, \Tdens \Slack = 0,\\
    &\int_{\Domain}\TPot \dx = 0.
  \end{aligned}
\end{equation}

Following the same steps described in~\cref{sec:discretization}, the (relaxed) discrete counterpart of equation in~\cref{eq:l2-shift} is a nonlinear system of equations with unknowns $(\Vect{\TPot},\Vect{\Tdens},\Vect{\Slack},\Vect{\Shift})\in (\REAL^{\Np}\times \REAL^{\Nr}\times \REAL^{\Nr}\times \REAL^{\Ntime})$ given by
\begin{alignat*}{2}
  \Fnewton_{\Vect{\Pot}}(\Vect{\TPot},\Vect{\Tdens})&=
  \left(
  \Block{1}{\Fnewton}_{\Vect{\Pot}};
  \ldots;
  \Block{\Ntimep}{\Fnewton}_{\Vect{\Pot}}
  \right)& =&\Vect{0}\in \REAL^{\Np}\ ,
  \\
  \New{\Fnewton}_{\Vect{\Tdens}}(\Vect{\TPot},\Vect{\Tdens},\Vect{\Slack},\Vect{\Shift})&=
  \left(
  \Block{1}{\New{\Fnewton}}_{\Vect{\Tdens}};
  \ldots;
  \Block{\Ntime}{\New{\Fnewton}}_{\Vect{\Tdens}}
  \right)& =&\Vect{0}\in \REAL^{\Nr}\ ,
  \\
  \Fnewton_{\Vect{\Slack}}(\Vect{\Tdens},\Vect{\Slack})&=
  \left(
  \Block{1}{\Fnewton}_{\Vect{\Slack}};
  \ldots;
  \Block{\Ntime}{\Fnewton}_{\Vect{\Slack}}
  \right)-\Relax\Vect{1}& =&\Vect{0}\in \REAL^{\Nr}\ ,
  \\
  \Fnewton_{\Vect{\Shift}}(\Vect{\TPot})&=
  \left(
  \Block{1}{\Fnewton}_{\Vect{\Shift}};
  \ldots;
  \Block{\Ntime}{\Fnewton}_{\Vect{\Shift}}
  \right)& =&\Vect{0}\in \REAL^{\Ntimep},
\end{alignat*}
where $\Fnewton^{\tstep}_{\Vect{\Pot}}$ and $\Fnewton^{\tstep}_{\Vect{\Slack}}$ are defined in~\cref{eq:conservation_discrete,eq:ortho_discrete}, while $\New{\Fnewton}^{\tstep}_{\Vect{\Tdens}}$ and $\Fnewton^{\tstep}_{\Vect{\Shift}}$ are given by
\begin{align*}
  \Block{\tstep}{\New{\Fnewton}}_{\Vect{\Tdens}}(
   \Vect{\TPot},
  \Vect{\Tdens},
  \Vect{\Slack},
  \Vect{\Shift})
  &:=
  \Block{\tstep}{\Fnewton}_{\Vect{\Tdens}}(
  \Vect{\Tdens},
  \Vect{\TPot},
  \Vect{\Slack})
  +
  \Vect{\AreaTdens} \Block{\tstep}{\Vect{\Shift}}\ ,\ &&\tstep=1,\ldots,\Ntime,
  \\
  \Block{\tstep}{\New{\Fnewton}}_{\Vect{\Shift}}(\Block{\tstep}{\Vect{\TPot}})
  &:=\Vect{\AreaPot}^{T}\Block{\tstep}{\Vect{\TPot}}\ ,\ &&\tstep=1,\ldots,\Ntimep,
\end{align*}
with $\Block{\tstep}{\Fnewton}_{\Vect{\Tdens}}$ as given in~\cref{eq:hj_discrete}. Note that there are $\Np+2\Nr+\Ntimep$ equations and $\Np+2\Nr+\Ntime$ unknowns. The additional equation fixes the global constant for the potential $\Vect{\TPot}$. The linear system arising from the Newton method becomes
\begin{align}
  \label{eq:newton-linsys-shift}
  \begin{pmatrix}
    \Ablock
    &
    \Bblock^T
    &
    \Matr{\ZeroMat}
    &
    \Matr{\ZeroMat}
    \\
    \Bblock
    &
    \Matr{\ZeroMat}
    &
    \MassTdensBlock
    &
    \MassTdensBlock\Matr{{\mathcal{\coarse{E}}}}^T
    \\
    \Matr{\ZeroMat}
    &
    \Diag{\Vect{\Slack}}
    &
    \Diag{\Vect{\Tdens}}
    &
    \Matr{\ZeroMat}
    \\
    \Matr{{\mathcal{E}}}  \Matr{\MassPotBlock} 
    &
    \Matr{\ZeroMat}
    &
    \Matr{\ZeroMat}
    &
    \Matr{\ZeroMat}
  \end{pmatrix}
  \begin{pmatrix}
    \Vect{\DTPot}
    \\
    \Vect{\DTdens}
    \\
    \Vect{\DSlack}
    \\
    \Vect{\DShift}
  \end{pmatrix}
  &
  =
  \begin{pmatrix}
    \Vect{f}
    \\
    \Vect{\New{g}}
    \\
    \Vect{h}
    \\
    \Vect{i}
  \end{pmatrix}
  =
  -
  \begin{pmatrix}
    \Vect{\Fnewton}_{\Vect{\Pot}}
    \\
    \New{\Vect{\Fnewton}}_{\Vect{\Tdens}}
    \\
    \Vect{\Fnewton}_{\Vect{\Slack}}-\Relax\Vect{1}
    \\
    \Vect{\Fnewton}_{\Vect{\Shift}}
  \end{pmatrix},
\end{align}
where matrices  $\Matr{\mathcal{{E}}}\in \REAL^{\Ntimep,(\Ntimep)\Npot}$  and  $\Matr{\mathcal{\coarse{E}}}\in \REAL^{\Ntime,\Ntime\Ntdens}$ are given by
\begin{equation*}
  \Matr{\mathcal{{E}}}=
  \begin{pmatrix}
    \Vect{1}_{{\Ncell}}^T & &
    \\
    & \ddots &
    \\
    &&\Vect{1}_{{\Ncell}}^T
  \end{pmatrix}
   , \quad
  \Matr{\mathcal{\coarse{E}}}=
  \begin{pmatrix}
    \Vect{1}_{\coarse{\Ncell}}^T & &
    \\
    & \ddots &
    \\
    &&\Vect{1}_{\coarse{\Ncell}}^T
  \end{pmatrix}
  .
\end{equation*}
Once the vectors $\Vect{\DTPot}$ and $\Vect{\DShift}$ in~\cref{eq:newton-linsys-shift} are found, $\Vect{\DPot}$ can be retrieved in the original system \cref{eq:saddle-point} using the formula
\begin{equation*}
  \label{eq:correction-discrete}
  \Block{\tstep}{\Vect{\DPot}}=\Block{\tstep}{\Vect{\DTPot}}
  +\sum_{i=1}^{\tstep}\Block{i}{\Deltat}\Block{i}{\Vect{\DShift}},
\end{equation*}
which is the finite-dimensional counterpart of~\cref{eq:pot-shift}.

\subsubsection{New saddle point linear system}
After eliminating the unknown $\Vect{\DSlack}$ as in~\cref{sec:linear-solver}, the linear system in~\cref{eq:newton-linsys-shift} becomes
\begin{align}
  \label{eq:newton-linsys-shift-reduced}
  \begin{pmatrix}
    \Ablock
    &
    \Bblock^T
    &
    \Matr{\ZeroMat}
    \\
    \Bblock
    &
    -\Cblock
    &
    \MassTdensBlock\Matr{{\mathcal{\coarse{E}}}}^T
    \\
    \Matr{{\mathcal{E}}}  \Matr{\MassPotBlock}
    &
    \Matr{\ZeroMat}
    &
    \Matr{\ZeroMat}
  \end{pmatrix}
  \begin{pmatrix}
    \Vect{\DTPot}
    \\
    \Vect{\DTdens}
    \\
    \Vect{\DShift}
  \end{pmatrix}
  &
  =
  \begin{pmatrix}
    \Vect{f}
    \\
    \Vect{{\tilde{g}}}
    \\
    \Vect{i}
  \end{pmatrix}.
\end{align}
The variable $\Vect{\DShift}$ can be eliminated as well by multiplying by $\Matr{\coarse{\Eblock}}$ the second block of equations in~\cref{eq:newton-linsys-shift}, yielding
\begin{equation*}
  \label{eq:remove-shift}
  \Matr{\coarse{\Eblock}} (\Bblock \Vect{\DTPot}-\Cblock\Vect{\DTdens})
  + \underbrace{\Matr{\coarse{\Eblock}}\MassTdensBlock\Matr{\coarse{\Eblock}}^T}_{|\Domain|\Identity{\Ntime}}\Vect{\DShift} = \Matr{\coarse{\Eblock}} \tilde{\Vect{g}}\,.
\end{equation*}
Plugging this expression in~\cref{eq:newton-linsys-shift-reduced}, the linear system in the unknowns $(\Vect{\DTPot};\Vect{\DTdens})$ becomes
\begin{equation}
  \label{eq:saddle-point-dtpot-drho}
  \begin{pmatrix}
    \Ablock
    &
    \Bblock^T
    \\
    \Matr{\ProjTdens} \Bblock
    &
    -\Matr{\ProjTdens} \Cblock
    \\
    \Matr{{\mathcal{E}}}  \Matr{\MassPotBlock}
    &
    \Matr{\ZeroMat}
  \end{pmatrix}
  \begin{pmatrix}
    \Vect{\DTPot}
    \\
    \Vect{\DTdens}
  \end{pmatrix}
  =
  \begin{pmatrix}
    \Vect{f}
    \\
    \Matr{\ProjTdens} \tilde{\Vect{g}}
    \\
    \Vect{i}
  \end{pmatrix}
  .
\end{equation}
where the matrix $\Matr{\ProjTdens}$, which is given by
\begin{equation*}
  \Matr{\ProjTdens}:=\Matr{\Identity{}}-\frac{1}{|\Domain|}\MassTdensBlock\Matr{{\mathcal{\coarse{E}}}}^T \Matr{\coarse{\Eblock}},
\end{equation*}
is a projector (indeed $\Matr{\ProjTdens}^2=\Matr{\ProjTdens}$). Moreover, its transpose $\Matr{\ProjTdens}^T$ is the discretization of the zero mean-projector, which maps a function $g:[0,1]\times\Domain\to\REAL$ to~$g-\int_{\Domain}g/|\Domain|$.

The linear system in~\eqref{eq:saddle-point-dtpot-drho} is not over-determined since $\Kernel{\Ablock}\subset\Kernel{\Matr{\ProjTdens}\Bblock}$ and the rows of $\Matr{{\mathcal{E}}}  \Matr{\MassPotBlock}$ are $\Ntimep$ vectors linearly independent from the rows of $\Ablock$ or $\Matr{\ProjTdens}\Bblock$.
Recalling~\cref{remark:zero-mean},
$
\Vect{\DTdens}\in \Kernel{\ProjTdens}
$
or equivalently
$
\Vect{\DTdens}\in \Image{\ProjTdens^T}
$
and, since $\ProjTdens^T$ is also a projector, we can write $\Vect{\DTdens} = \ProjTdens^T \Vect{\DTdens}$. Thus, we can rewrite system \eqref{eq:saddle-point-dtpot-drho} as
\begin{align}
  \label{eq:newton-linsys-shift-reduced-projected}
  \begin{pmatrix}
    \Ablock
    &
    \Bblock^T\Matr{\ProjTdens}^T
    \\
    \Matr{\ProjTdens} \Bblock
    &
    -\Matr{\ProjTdens} \Cblock\Matr{\ProjTdens}^T
  \end{pmatrix}
  \begin{pmatrix}
    \Vect{\delta\hat{\Pot}}
    \\
    \Vect{\DTdens}
  \end{pmatrix}
  &
  =
  \begin{pmatrix}
    \Vect{f}
    \\
    \Matr{\ProjTdens}\Vect{{\tilde{g}}}
  \end{pmatrix},
\end{align}
where we neglect the constraints $\Matr{{\mathcal{E}}}  \Matr{\MassPotBlock}\Vect{\DTPot}=\Vect{i}$. For any solution $( \Vect{\delta\hat{\Pot}};\Vect{\delta {\Tdens}})$, the solution $(\Vect{{\DTPot}};\Vect{{\DTdens}})$ of~\cref{eq:saddle-point-dtpot-drho} can be recovered by simply correcting the mean of $\Vect{\delta\hat{\Pot}}$.

The linear system in~\cref{eq:newton-linsys-shift-reduced-projected} does not suffer from the first issue mentioned above associated to the singularity of matrix $\Ablock$. In fact, note that $\Image{\Bblock^T \Matr{\ProjTdens}^T}\perp \Kernel{\Ablock}$, since $\Kernel{\Ablock}\subset\Kernel{\Matr{\ProjTdens}\Bblock}$. This ensures that the dual Schur complement $\Matr{\SchurDual}=-\Matr{\ProjTdens}(\Cblock + \Bblock \Pinv{\Ablock} \Bblock^T) \Matr{\ProjTdens}^T$ (where $\Pinv{\Ablock}$ denotes the Moore-Penrose pseudo-inverse of $\Ablock$) is well defined. Moreover the following (ideal) block triangular preconditioner  
\begin{gather}
  \label{eq:triang}
  \Matr{P}_{t}
  =
  \begin{pmatrix}
    \Matr{\Amat}& \Bblock^T \Matr{\ProjTdens}^T
    \\
    \ZeroMat&  -\Matr{\SchurDual}
  \end{pmatrix}
\end{gather}
is well defined as well. In fact, given a vector $(\Vect{\RhsPot};\Vect{\RhsTdens})$ with $\Vect{\RhsPot}\in \Kernel{\Ablock}^\perp$ and $\Vect{\RhsTdens}\in\Image{\ProjTdens}$, the application of this preconditioner requires to compute the vector $(\Vect{x};\Vect{y})$ solving the following two linear systems
\begin{align}
  \label{eq:wls}
  -\Matr{\SchurDual} \Vect{y} &= \Vect{\RhsTdens}\,,
  \\
  \label{eq:x-triang}
  \Ablock \Vect{x} & =  \Vect{b}=\Vect{\RhsPot} - \Bblock^T \ProjTdens^T \Vect{y}\,.
\end{align}
The linear system in~\cref{eq:x-triang} is now well defined, since   $\Vect{b}\in\Kernel{\Ablock}^\perp$ at any preconditioner application. In fact, both vectors $\Bblock^T \ProjTdens^T \Vect{y}$ and $\Vect{\RhsPot}$ belong to $\Kernel{\Ablock}^\perp$, the first vector because $\Image{\Bblock^T  \Matr{\ProjTdens}^T}\subset\Kernel{\Ablock}^\perp$, the second vector since it is a linear combination of the vector $\Vect{f}$ in~\cref{eq:newton-linsys-shift-reduced-projected} and vectors in $\Image{\Ablock}\cup\Image{\Bblock ^T\ProjTdens^T}$, all belonging to $\Kernel{\Ablock}^\perp$.

\subsubsection{Approximating the inverse of the dual \Schur{} complement}
\label{sec:bb_definition}
We turn now to the issue of approximating the Schur complement. Our first attempt was to approximate it with $\Matr{\ProjTdens}(\Cblock + \Bblock\Inv{\Diag{\Ablock}} \Bblock^T) \Matr{\ProjTdens}^T$, but it led to poor performance. We tried to apply the Algebraically stabilized Least Square Commutator proposed in~\cite[Section 4.2]{Elman2007}, but this leads to poor results, for different combinations of the relaxation parameters that appear in such approach. We attribute this phenomenon to the fact that the matrix $\Cblock$ is not a stabilization operator for the saddle point system, as assumed in~\cite{Elman2007}.

In order to devise a better approximation of the \Schur{} complement, here we follow the idea in~\cite{Elman2006,Elman2007,elman2005preconditioning} of looking at the infinite-dimensional differential operators associated to the matrices that compose it and try to deduce a possible approximation of its inverse. The following proposition shows a differential identity that goes in this direction.
\begin{Prop}
  \label{prop:space-commute}
  Let $\Tdens:[0,1]\times\Domain\to\REALP$, $\Pot:[0,1]\times\Domain\to\REAL$, smooth enough.
  Denoting by $-\Delta_{\Tdens}=-\Div(\Tdens \Grad)$, the following operator identity holds:
     \begin{equation}
      \label{eq:commuting}
      \begin{split}
        \left( \Dt{}+\Div(\cdot \Grad \Pot) \right) (-\Delta_{\Tdens})
        =
        - \Delta_{\Tdens}
        \left(\Dt +\Grad \Pot \cdot \Grad \right)
        -\Div
        \left(
        \left(\Dt{\Tdens} + \Div(\Tdens \Grad \Pot)
        \right)
        \Grad
        \right)
        \\
        +2\Div(\Tdens \Grad^2 \Pot \Grad )\,. 
      \end{split}
    \end{equation}
\end{Prop}
\begin{proof}
  Consider any function $g_0:[0,1]\times\Domain \to \REAL$, such that $\int_{\Domain} g_0(t,\cdot)=0$ for all $t\in[0,1]$ and let us define $\Sol:[0,1]\times\Domain \to \REAL$ given by
    \begin{gather}
    \label{eq:div} 
    -\Div(\Tdens \Grad \Sol) = g_0
    \end{gather}
    with zero Neumann boundary conditions.
  Multiplying both sides of equation~\cref{eq:div} by $\partial_{i}\Pot$ (where $\partial_i=\partial_{x_i}$) and taking the divergence, we obtain 
  \begin{equation*}
    -\sum_{i,j} \partial_{i} \left(
    \partial_{i} \Pot \partial_{j} (\Tdens \partial_{j}\Sol)
    \right)
    =
    \Div(g_0\Grad \Pot).
  \end{equation*}
  Using the identity $
  \partial_{i} \Pot\partial_{j} (\Tdens \partial_{j}\Sol) = \partial_{j}(\partial_{i} \Pot \Tdens \partial_{j} \Sol) - \partial_{j}(\partial_{i} \Pot)\Tdens \partial_{j}\Sol,
  $
  we get
   \begin{align*}
     \Div(g_0\Grad \Pot)
     &=
     -
     \sum_{i,j} \partial_{i} \left(
     \partial_{j} (\partial_{i} \Pot \Tdens \partial_{j}\Sol)
     -  
     \partial_{j,i} \Pot \Tdens \partial_{j}\Sol
     \right)
     \\
     &=
      -
    \sum_{j,i} \partial_{j} \left(
    \partial_{i,j} \Sol \Tdens \partial_{i}\Pot  
    \right) 
     -
    \sum_{j,i}        
    \partial_{j} \left(
    \partial_{i}
    (\Tdens \partial_{i}\Pot)
    \partial_{j} \Sol 
    \right)
    +\Div(  \Tdens \Grad^2 \Pot \Grad\Sol)
    \\    
    &=
    \!-\!
    \sum_{j} \partial_{j} \left(
    \sum_{i} \partial_{i,j} \Sol \Tdens \partial_{i}\Pot
    \right)
    \!-\!
    \sum_{j}        
    \partial_{j} \left(
    \sum_{i}
    \partial_{i}
    (\Tdens \partial_{i}\Pot)
    \partial_{j} \Sol
    \right)
    \\
    & \qquad \qquad \qquad +\Div(  \Tdens \Grad^2 \Pot \Grad \Sol)
    \\
     &=
     -\Div(\Tdens \Grad^2\Sol \Grad \Pot)
     -\Div(\Div(\Tdens \Grad \Pot)    \Grad \Sol)
     +\Div( \Tdens \Grad^2 \Pot\Grad\Sol)\,.
  \end{align*}
   Using the identity $\Grad (\Grad \Pot \cdot \Grad \Sol)=\Grad^2 \Pot\Grad \Sol+\Grad^2\Sol\Grad\Pot$, we obtain
   \begin{equation}
     \label{eq:commute-space}
     \Div(g_0\Grad \Pot)
     =
     -\Div(\Tdens \Grad (\Grad \Pot\cdot \Grad \Sol))
     -\Div(\Div(\Tdens \Grad \Pot)    \Grad \Sol)     
     +2\Div( \Tdens \Grad^2 \Pot\Grad\Sol)\,.
   \end{equation}
   Now, differentiating with respect to time the expression in~\cref{eq:div}, we obtain 
   \begin{equation}
  \label{eq:commute-time}
  -\Div\left(\Dt{\Tdens} \Grad \Sol + \Tdens \Grad \Dt{\Sol} \right) =\Dt{g_0}\,.
   \end{equation}
   Combining~\cref{eq:commute-space}, \cref{eq:commute-time} we obtain
   \begin{equation*}
   \begin{split}
     \Dt{g_0}+\Div(g_0\Grad \Pot)
     =
     -\Div(\Tdens \Grad (\Dt{\Sol} + \Grad \Pot\cdot \Grad \Sol))
     -\Div\left(\left(\Dt{\Tdens} + \Div(\Tdens \Grad \Pot) \right) \Grad \Sol\right)
     \\
     +2\Div( \Tdens \Grad^2 \Pot\Grad\Sol)\,.
   \end{split}
   \end{equation*}
   Since $g_0=-\Delta_\Tdens \Sol$ by~\cref{eq:div}, we proved ~\cref{eq:commuting}.
\end{proof}

Our approximation of the Schur complement is based on neglecting the terms
$-\Div\left(\left(\Dt{\Tdens} +\Div(\Tdens \Grad \Pot\right)\right) \Grad$  and $2\Div(\Tdens \Grad^2 \Pot \Grad)$ on the right-hand side of~\cref{eq:commuting}.
The first term may be assumed to be small since it contains the continuity equation, which is the first nonlinear equation in the system~\cref{eq:l2-pde}. This term is small at the initial and final Newton steps within each \IP{} iteration, and experimental observations showed it remains small during the intermediate Newton iterations. The second term can be neglected under certain assumptions. If the optimal transport is a translation this term is null. When the transported measures $\Source$ and $\Sink$ are Gaussian \cite{Hessian2007}, or more generally log-concave probability densities \cite{caffarelli2000monotonicity,colombo2017lipschitz}, the largest eigenvalue of $\Grad^2 \Pot(0,x)$ is bounded. However, note that in general this term may be non null, as for example in the test case 3. {Supposing that both terms can be neglected, we get the following approximate identity:
\begin{equation*}
  \label{eq:commute}
  \left(\Dt{} + \Div (\cdot \Grad \Pot)\right)(-\Delta_{\Tdens}) \approx  -\Delta_{\Tdens} \left( \Dt{} + \Grad \Pot \cdot \Grad \right).
\end{equation*}
In our finite-dimensional problem, this expression translates into the approximate matrix identity
\begin{equation}
  \label{eq:commute-matrix}
  - \Bblock^T \MassTdensBlock^{-1}  \AblockTdens \approx \Ablock \MassPotBlock^{-1} \tildeBblock\,,
\end{equation}
where $\AblockTdens\in\REAL^{\Nr,\Nr}$ and $\tilde{\Bblock}\in\REAL^{\Np,\Nr}$ discretize the operators $-\Delta_\Tdens$ and $\Dt{}+\Grad \Pot \cdot \Grad$, respectively, similarly to the matrices $\Ablock$ and $\Bblock$ but with different dimensions. The scaling factors $\MassTdensBlock^{-1}$ and $\MassPotBlock$ in~\cref{eq:commute-matrix} are introduced to match the scaling dimensions of the matrices $\Ablock$, $\Bblock$, and  $\Bblock^T$.} The most natural definition of $\AblockTdens$ and $\tilde{\Bblock}$ is the following:
\begin{gather*}
  \AblockTdens=\BlkDiag{
  \left(\Block{\tstep}{\tilde{\Matr{A}}}\right)_{\tstep=1}^{\Ntime}
  }
  , \quad 
  \Block{\tstep}{\tilde{\Matr{A}}}:=\Matr{\Mdiv{\coarse{\Triang},\coarse{\Edges}}}  \Diag{\Recon{\coarse{\Edges}}\Of{\Block{\tstep}{\Vect{\Tdens}}}} \Matr{\Mgrad{\coarse{\Edges},\coarse{\Triang}}},
  \\
  \tilde{\Bblock}:= -\MassPotBlock \TdensToPotBlock{}\tilde{\Matr{\TimeDiff}}^T + \Matr{\Gpot} \Matr{\tilde{\Gradblock}}  \TdensToPotBlock{}^T \TimeAvgBlock^T,
\end{gather*}
  where the matrices $\TdensToPotBlock{}$, $\Matr{\TimeAvgBlock}$, and $\Matr{\Gpot}$ are defined in~\cref{eq:dt-it,eq:timeavg-G}, while the matrices $\Matr{\tilde{\TimeDiff}}\in \REAL^{\Ntime\Ntdens,(\Ntimep)\Ntdens}$ and $\tilde{\Matr{\Gradblock}}\in \REAL^{\Ntime\coarse{\Nedge},\Ntime\Ntdens}$ are
\[
\Matr{\tilde{\TimeDiff}}=\frac{1}{\Deltat}
 \begin{pmatrix}
    -\Identity{\Ntdens}
    &
    \Identity{\Ntdens}
    &
    \ZeroMat
    &
    \ZeroMat
    \\
    \ZeroMat
    &
    \ddots
    &
    \ddots
    &
    \ZeroMat
    \\
    \ZeroMat
    &
    \ZeroMat
    &
    -\Identity{\Ntdens}
    &
    \Identity{\Ntdens}
 \end{pmatrix}
 , 
  \quad
\tilde{\Matr{\Gradblock}}=\BlkDiag{(\Mgrad{\coarse{\Edges},\coarse{\Triang}})_{\tstep=1}^{\Ntime}}.
\]

Using the approximate identity in~\cref{eq:commute-matrix}, the linear system in~\cref{eq:wls} becomes
\begin{equation}
  \label{eq:schur-linear-system}
  -\Vect{d}=\Sblock \Vect{y}=\ProjTdens\left(\Cblock + \Bblock\Pinv{\Ablock}\Bblock^T\right)\ProjTdens^T\Vect{y}
  \approx\ProjTdens(\Cblock - \Bblock\MassPotBlock^{-1} \Tilde{\Bblock}
  \Pinv{\AblockTdens} \MassTdensBlock)\ProjTdens^T \Vect{y}\,.
\end{equation}
Since $\Kernel{\ProjTdens  \MassTdensBlock}=\Image{\Eblock^T}$ are discrete functions constant in space, we have that $ \Image{\MassTdensBlock \ProjTdens^T} \subset \Image{\AblockTdens}$ and 
the last expression is equal to
\begin{align*}
  -\Vect{d}=\ProjTdens(\Cblock \MassTdensBlock^{-1}
  \AblockTdens
  \Pinv{\AblockTdens} \MassTdensBlock
  -
  \Bblock\MassPotBlock^{-1}
  \Tilde{\Bblock}\Pinv{\AblockTdens}\MassTdensBlock)\ProjTdens^T\Vect{y}
  \,.
\end{align*}
Hence, we only have to solve the linear system
\begin{equation*}
-\Vect{d}=\ProjTdens(\Cblock\MassTdensBlock^{-1}\AblockTdens - \Bblock\MassPotBlock^{-1}\Tilde{\Bblock}) \Pinv{\AblockTdens } \MassTdensBlock  \ProjTdens^T\Vect{y}\,,
\end{equation*}
which is compatible since $\Vect{d}\in \Image{\ProjTdens}$. Note that the dense matrix  $(\Cblock + \Bblock\Pinv{\Ablock}\Bblock^T)$ in~\cref{eq:schur-linear-system} is replaced by the operator $(\Cblock\MassTdensBlock^{-1}\AblockTdens - \Bblock\MassPotBlock^{-1}\Tilde{\Bblock}) \Pinv{\AblockTdens } \MassTdensBlock$ which is dense by the presence of the pseudo-inverse of $\AblockTdens$. But now this operator is composed by two factors. This means that we can approximate the action of their respective inverse by solving a sparse linear system followed by a matrix-vector product. In fact, since the solution $\Vect{\DTdens}$ in~\cref{eq:newton-linsys-shift-reduced-projected} belongs to the image of $\ProjTdens^T$, we may look for a solution $\Vect{y}\in \Image{ \ProjTdens^T}$. By setting
\[\Vect{z}=\Pinv{\AblockTdens } \MassTdensBlock  \ProjTdens^T\Vect{y}=
\Pinv{\AblockTdens } \MassTdensBlock\Vect{y}\,,
\]
we can compute $\Vect{y}$ by first solving the linear system
\begin{equation}
  \label{eq:pseudoschur}
  -\Vect{d}=(\Cblock\MassTdensBlock^{-1}\AblockTdens - \Bblock\MassPotBlock^{-1}\Tilde{\Bblock}) \Vect{z}
\end{equation}
and then computing $\Vect{y}=\MassTdensBlock^{-1} \AblockTdens\Vect{z}$.

Moreover, in order to cope with the ill-conditioning derived from the presence of the matrix $\Cblock=\MassTdensBlock\Diag{\Vect{\Slack}/\Vect{\Tdens}}$, we scale~\cref{eq:pseudoschur} by $\Diag{\Vect{\Tdens}}$ and we solve the linear system 
\begin{equation}
  \label{eq:pseudoschur-bis}
  (\Diag{\Vect{\Slack}} \AblockTdens - \Diag{\Vect{\Tdens}}\Bblock\MassPotBlock^{-1}\Tilde{\Bblock}) \Vect{z}=
  -\Diag{\Vect{\Tdens}}\Vect{d},
\end{equation}
similarly to what done in~\cref{eq:linsys-simple-mode2} for the SIMPLE preconditioner.

We refer to the resulting preconditioner as ``$\BB$-preconditioner'', due to presence of this composed operator in the approximate factorization of the Schur complement. This choice was done in analogy with the ``$\Matr{B}\Matr{F}\Matr{B}t$-preconditioner'' introduced in~\cite{elman2005preconditioning}, whose ideas inspired the preconditioner described in this section.

\subsubsection{Numerical results}
In this section we present the results obtained on test case 2 using the preconditioner described in~\cref{sec:bb_definition}. Its application requires the approximate solution of the linear systems in~\cref{eq:pseudoschur-bis,eq:x-triang}. The first is solved using the AGMG solver with inner tolerance $\InnerTol^{\Sblock}=1e-1$ (determined experimentally). The second involves the solution of $\Ntime+1$ weighted Laplacian systems. We solved them with the AGMG solver with tolerance $\InnerTol^{\Ablock}=5e-2$. This value, slightly smaller than those adopted in previous cases, improves the robustness of the preconditioner. The number of multigrid iterations for these systems ranges between 2 and 3 in all cases.

\addtolength{\tabcolsep}{-2.5pt}
\begin{table}
  \begin{center}
    {\footnotesize
      

\begin {tabular}{R{0.8cm }|R{0.44cm }R{0.44cm }R{0.44cm }R{0.5cm }|R{0.68cm }R{0.68cm }R{0.82cm }R{0.82cm }|R{0.62cm }R{0.62cm }R{0.62cm }R{0.62cm }}
\toprule
 
    $\Ntime+1$ 
    & 16 & 32 & 64 & 128  
    & 16 & 32 & 64 & 128
    & 16 & 32 & 64 & 128 
    \\
    \midrule
    
    $\Relax$
    & \multicolumn {4}{c|}{\textbf {\AvgOuter}}
    & \multicolumn {4}{c|}{\textbf {\AvgCPU }}
    & \multicolumn {4}{c}{\textbf {\AvgInner }}
    \\
    \midrule
     
    \multicolumn {13}{c}{$\Triang^{0} (
    \Ncell=224,
    \coarse{\Ncell}=672)$ \quad \#DOF=(1.4e4, 2.8e4, 5.7e4, 1.1e5)}
    \\
    \midrule
    1  & \cellcolor {gray!20}7 & 9 & 10 & 12 & \cellcolor {gray!20}1.2e0 & 2.5e0 & 5.4e0 & 1.4e1 & \cellcolor {gray!20}2.1 & 2.3 & 2.6 & 3.2\\ 
2e-1 & \cellcolor {gray!20}8 & 9 & 9 & 10 & \cellcolor {gray!20}1.3e0 & 2.6e0 & 5.5e0 & 1.3e1 & \cellcolor {gray!20}2.6 & 3.2 & 3.9 & 4.8\\ 
4e-2 & \cellcolor {gray!20}12 & 13 & 14 & 15 & \cellcolor {gray!20}1.3e0 & 2.8e0 & 5.7e0 & 1.4e1 & \cellcolor {gray!20}2.8 & 3.6 & 4.5 & 5.5\\ 
8e-3 & \cellcolor {gray!20}17 & 18 & 21 & 25 & \cellcolor {gray!20}1.4e0 & 3.0e0 & 6.4e0 & 1.8e1 & \cellcolor {gray!20}3.1 & 4.1 & 5.0 & 6.3\\ 
2e-3 & \cellcolor {gray!20}20 & 25 & 33 & 44 & \cellcolor {gray!20}1.5e0 & 3.3e0 & 7.5e0 & 2.3e1 & \cellcolor {gray!20}3.6 & 4.4 & 5.1 & 6.2\\ 
3e-4 & \cellcolor {gray!20}30 & 42 & 63 & 88 & \cellcolor {gray!20}1.6e0 & 3.9e0 & 1.0e1 & 3.4e1 & \cellcolor {gray!20}3.7 & 4.5 & 4.7 & 5.4\\ 
6e-5 & \cellcolor {gray!20}46 & 72 & 117 & 180 & \cellcolor {gray!20}1.9e0 & 4.9e0 & 1.5e1 & 5.9e1 & \cellcolor {gray!20}3.9 & 4.1 & 4.4 & 4.4\\ 
1e-5 & \cellcolor {gray!20}76 & 127 & 248 & \dag & \cellcolor {gray!20}2.4e0 & 6.8e0 & 2.8e1 & \dag & \cellcolor {gray!20}4.1 & 4.3 & 4.0 & \dag\\ 
3e-6 & \cellcolor {gray!20}136 & 324 & \dag & \dag & \cellcolor {gray!20}3.5e0 & 1.5e1 & \dag & \dag & \cellcolor {gray!20}4.1 & 4.8 & \dag & \dag\\ 
5e-7 & \cellcolor {gray!20}297 & \dag & \dag & \dag & \cellcolor {gray!20}6.2e0 & \dag & \dag & \dag & \cellcolor {gray!20}4.3 & \dag & \dag & \dag\\ 
\toprule  & \cellcolor {gray!20}58 & \dag & \dag & \dag & \cellcolor {gray!20}2.1e0 & \dag & \dag & \dag & \cellcolor {gray!20}4.0 & \dag & \dag & \dag\\ 
\bottomrule 
 \end{tabular}
      

\begin {tabular}{R{0.8cm }|R{0.44cm }R{0.44cm }R{0.44cm }R{0.5cm }|R{0.68cm }R{0.68cm }R{0.82cm }R{0.82cm }|R{0.62cm }R{0.62cm }R{0.62cm }R{0.62cm }}
\toprule
 
    \multicolumn {13}{c}{$\Triang^{1} (
    \Ncell=896,
    \coarse{\Ncell}=2688)$ \quad \#DOF=(5.6e4, 1.1e5, 2.3e5, 4.6e5)}
    \\
    \midrule
    1  & 7 & \cellcolor {gray!20}8 & 9 & 10 & 1.7e0 & \cellcolor {gray!20}3.4e0 & 8.1e0 & 2.1e1 & 2.9 & \cellcolor {gray!20}2.3 & 2.5 & 3.2\\ 
2e-1 & 8 & \cellcolor {gray!20}8 & 9 & 10 & 1.6e0 & \cellcolor {gray!20}3.4e0 & 8.0e0 & 2.1e1 & 2.4 & \cellcolor {gray!20}2.8 & 3.6 & 4.3\\ 
4e-2 & 13 & \cellcolor {gray!20}15 & 15 & 15 & 1.8e0 & \cellcolor {gray!20}4.0e0 & 1.0e1 & 2.6e1 & 2.3 & \cellcolor {gray!20}3.3 & 4.2 & 5.4\\ 
8e-3 & 17 & \cellcolor {gray!20}20 & 21 & 21 & 1.9e0 & \cellcolor {gray!20}4.7e0 & 1.2e1 & 3.3e1 & 2.6 & \cellcolor {gray!20}3.6 & 4.8 & 5.8\\ 
2e-3 & 22 & \cellcolor {gray!20}26 & 27 & 30 & 2.2e0 & \cellcolor {gray!20}5.5e0 & 1.5e1 & 4.4e1 & 2.9 & \cellcolor {gray!20}4.1 & 5.4 & 5.8\\ 
3e-4 & 29 & \cellcolor {gray!20}34 & 38 & 44 & 2.6e0 & \cellcolor {gray!20}6.5e0 & 2.0e1 & 5.9e1 & 3.2 & \cellcolor {gray!20}4.4 & 5.7 & 6.0\\ 
6e-5 & 37 & \cellcolor {gray!20}46 & 61 & 86 & 3.1e0 & \cellcolor {gray!20}8.7e0 & 2.8e1 & 1.1e2 & 3.7 & \cellcolor {gray!20}4.7 & 5.8 & 6.3\\ 
1e-5 & 52 & \cellcolor {gray!20}65 & 88 & 157 & 4.0e0 & \cellcolor {gray!20}1.3e1 & 4.3e1 & 2.0e2 & 4.2 & \cellcolor {gray!20}5.4 & 8.1 & 7.7\\ 
3e-6 & 72 & \cellcolor {gray!20}100 & 188 & 395 & 5.0e0 & \cellcolor {gray!20}1.7e1 & 1.0e2 & 4.9e2 & 4.5 & \cellcolor {gray!20}5.8 & 11.0 & 8.1\\ 
5e-7 & 100 & \cellcolor {gray!20}154 & \dag & \dag & 6.7e0 & \cellcolor {gray!20}2.5e1 & \dag & \dag & 4.9 & \cellcolor {gray!20}6.0 & \dag & \dag\\ 
\toprule  & 32 & \cellcolor {gray!20}43 & \dag & \dag & 2.9e0 & \cellcolor {gray!20}8.4e0 & \dag & \dag & 4.0 & \cellcolor {gray!20}5.1 & \dag & \dag\\ 
\bottomrule 
 \end{tabular}
      

\begin {tabular}{R{0.8cm }|R{0.44cm }R{0.44cm }R{0.44cm }R{0.5cm }|R{0.68cm }R{0.68cm }R{0.82cm }R{0.82cm }|R{0.62cm }R{0.62cm }R{0.62cm }R{0.62cm }}
\toprule
 
    \multicolumn {13}{c}{$\Triang^{2} (
    \Ncell=3584,
    \coarse{\Ncell}=10752)$ \quad \#DOF=(2.3e5, 4.6e5, 9.1e5, 1.8e6)}
    \\
    \midrule
    1  & 6 & 7 & \cellcolor {gray!20}8 & 10 & 2.7e0 & 6.6e0 & \cellcolor {gray!20}1.7e1 & 5.6e1 & 3.5 & 2.9 & \cellcolor {gray!20}2.4 & 2.6\\ 
2e-1 & 9 & 9 & \cellcolor {gray!20}10 & 10 & 3.0e0 & 7.3e0 & \cellcolor {gray!20}1.9e1 & 5.9e1 & 2.6 & 2.5 & \cellcolor {gray!20}2.9 & 3.6\\ 
4e-2 & 14 & 16 & \cellcolor {gray!20}18 & 17 & 3.9e0 & 1.0e1 & \cellcolor {gray!20}2.9e1 & 8.6e1 & 2.9 & 2.6 & \cellcolor {gray!20}3.5 & 4.4\\ 
8e-3 & 19 & 22 & \cellcolor {gray!20}25 & 25 & 4.9e0 & 1.3e1 & \cellcolor {gray!20}4.0e1 & 1.2e2 & 3.3 & 2.7 & \cellcolor {gray!20}3.9 & 4.9\\ 
2e-3 & 25 & 28 & \cellcolor {gray!20}32 & 30 & 7.0e0 & 1.6e1 & \cellcolor {gray!20}5.0e1 & 1.5e2 & 3.8 & 2.9 & \cellcolor {gray!20}4.3 & 5.7\\ 
3e-4 & 32 & 34 & \cellcolor {gray!20}42 & 41 & 8.7e0 & 2.0e1 & \cellcolor {gray!20}6.5e1 & 2.0e2 & 4.4 & 3.2 & \cellcolor {gray!20}4.3 & 5.8\\ 
6e-5 & 44 & 43 & \cellcolor {gray!20}53 & 58 & 1.2e1 & 2.5e1 & \cellcolor {gray!20}8.0e1 & 2.8e2 & 4.8 & 3.4 & \cellcolor {gray!20}4.4 & 6.2\\ 
1e-5 & 62 & 61 & \cellcolor {gray!20}84 & 103 & 1.6e1 & 3.6e1 & \cellcolor {gray!20}1.3e2 & 4.9e2 & 5.4 & 3.9 & \cellcolor {gray!20}4.7 & 6.4\\ 
3e-6 & 92 & 94 & \cellcolor {gray!20}130 & 186 & 2.5e1 & 5.8e1 & \cellcolor {gray!20}2.0e2 & 9.6e2 & 6.7 & 4.3 & \cellcolor {gray!20}5.5 & 8.8\\ 
5e-7 & 125 & 143 & \cellcolor {gray!20}211 & \dag & 3.6e1 & 9.5e1 & \cellcolor {gray!20}3.5e2 & \dag & 7.3 & 5.2 & \cellcolor {gray!20}6.5 & \dag\\ 
\toprule  & 35 & 41 & \cellcolor {gray!20}54 & \dag & 9.7e0 & 2.6e1 & \cellcolor {gray!20}8.7e1 & \dag & 5.5 & 4.0 & \cellcolor {gray!20}5.2 & \dag\\ 
\bottomrule 
 \end{tabular}
      

\begin {tabular}{R{0.8cm }|R{0.44cm }R{0.44cm }R{0.44cm }R{0.5cm }|R{0.68cm }R{0.68cm }R{0.82cm }R{0.82cm }|R{0.62cm }R{0.62cm }R{0.62cm }R{0.62cm }}
\toprule
 
    \multicolumn {13}{c}{$\Triang^{3} (
    \Ncell=14336,
    \coarse{\Ncell}=43008)$ \quad \#DOF=(9.0e5, 1.8e6, 3.7e6, 7.3e6)}
    \\
    \midrule
    1  & 7 & 7 & 8 & \cellcolor {gray!20}9 & 8.3e0 & 2.1e1 & 5.7e1 & \cellcolor {gray!20}2.2e2 & 4.3 & 3.6 & 2.9 & \cellcolor {gray!20}2.4\\ 
2e-1 & 9 & 9 & 10 & \cellcolor {gray!20}11 & 9.8e0 & 2.5e1 & 7.0e1 & \cellcolor {gray!20}2.6e2 & 3.1 & 2.8 & 2.7 & \cellcolor {gray!20}2.8\\ 
4e-2 & 16 & 17 & 18 & \cellcolor {gray!20}20 & 1.5e1 & 4.0e1 & 1.1e2 & \cellcolor {gray!20}4.4e2 & 3.3 & 3.2 & 2.9 & \cellcolor {gray!20}3.7\\ 
8e-3 & 25 & 29 & 30 & \cellcolor {gray!20}32 & 2.7e1 & 6.3e1 & 1.7e2 & \cellcolor {gray!20}7.3e2 & 4.2 & 3.4 & 2.7 & \cellcolor {gray!20}3.8\\ 
2e-3 & 34 & 36 & 39 & \cellcolor {gray!20}48 & 5.0e1 & 8.1e1 & 2.3e2 & \cellcolor {gray!20}1.1e3 & 4.8 & 3.7 & 2.9 & \cellcolor {gray!20}4.2\\ 
3e-4 & 47 & 43 & 52 & \cellcolor {gray!20}58 & 7.6e1 & 1.0e2 & 3.0e2 & \cellcolor {gray!20}1.3e3 & 4.9 & 4.0 & 2.9 & \cellcolor {gray!20}4.3\\ 
6e-5 & 70 & 54 & 59 & \cellcolor {gray!20}73 & 1.2e2 & 1.4e2 & 3.4e2 & \cellcolor {gray!20}1.6e3 & 5.1 & 4.7 & 3.2 & \cellcolor {gray!20}4.5\\ 
1e-5 & 110 & 68 & 72 & \cellcolor {gray!20}100 & 1.9e2 & 2.0e2 & 4.3e2 & \cellcolor {gray!20}2.3e3 & 5.6 & 5.6 & 3.7 & \cellcolor {gray!20}4.6\\ 
3e-6 & 146 & 105 & 107 & \cellcolor {gray!20}149 & 2.8e2 & 3.2e2 & 6.5e2 & \cellcolor {gray!20}3.4e3 & 6.4 & 6.4 & 4.1 & \cellcolor {gray!20}4.9\\ 
5e-7 & 190 & 162 & 161 & \cellcolor {gray!20}237 & 3.9e2 & 6.1e2 & 1.0e3 & \cellcolor {gray!20}5.6e3 & 7.3 & 8.5 & 5.1 & \cellcolor {gray!20}5.5\\ 
\toprule  & 55 & 44 & 52 & \cellcolor {gray!20}68 & 9.4e1 & 1.3e2 & 3.2e2 & \cellcolor {gray!20}1.6e3 & 5.8 & 5.7 & 3.9 & \cellcolor {gray!20}4.8\\ 
\bottomrule 
 \end{tabular}
    }
  \end{center}
  \caption{
    Numerical results using the $\BB$-preconditioner. Each sub-table reports the results for a given mesh, from the coarsest (top) to the finest (bottom), and for different time discretizations ($\Deltat=1/(\Ntime+1)$). Each sub-table reports, while $\Relax$ is reduced (leftmost column), the averaged outer iteration and CPU time per linear system (\AvgOuter{} and \AvgCPU{}) and the averaged inner iterations per outer iteration \AvgInner{} for solving the linear system in~\cref{eq:primal-linsys} (metrics defined in~\cref{eq:pseudoschur}). A final row summarizes the averages on the whole simulation. We highlighted in gray the time-space combination providing the uniform discretization. The $\dag$ symbol denotes those \IP{} steps where the linear solver failed.
     }
  \label{tab:sol20}
\end{table}
\addtolength{\tabcolsep}{2.5pt}
In~\cref{tab:sol20} we summarized the results for the metrics defined in~\cref{sec:preconditioning} that we used to measure the preconditioner performance. The number of averaged inner iterations \AvgInner{} refers to the solution of the linear system~\cref{eq:pseudoschur-bis}. This number remains bounded between 2.1 and 8.5 iterations per preconditioner applications, it tends to increase slightly with $\Relax \to 0$ but it is not affected by the size of the problem.

Unfortunately, the number of outer iterations \AvgOuter{} tends to increase as $\mu\to 0$, exceeding in some cases the $400$ limit we fixed for the FGMRES iterations. This phenomenon is more pronounced when the temporal scale is finer than the spatial one to. This loss in efficiency for $\mu\to 0$ is more evident in the test cases $2$ and $3$ whose solutions have compact support, while for the test case $1$ it is less pronounced. This suggests that this phenomenon can be related to the degeneracy of the linear system~\cref{eq:pseudoschur-bis}, which becomes underdetermined due to the presence of terms that go to zero on both sides of the equation.

Nevertheless, for $\Relax\approx 1e-5$,  this preconditioner is the only one scaling well with respect to the temporal and spatial discretization. The number of averaged outer iterations \AvgOuter{} increases only slightly when the mesh is refined or the time steps is halved.  This results in a CPU time that scales sightly worse than linearly with respect to the number of degrees of freedom used to discretize the problem.

\subsection{Summary of numerical results}
\label{sec:resume}
In this section, we summarize the numerical results to compare the pros and cons of the preconditioning approaches proposed in this paper.
\begin{figure}
  \centerline{
    \includegraphics[width=0.8\textwidth,
      trim={0.0cm 0.0cm 0cm 0.0cm},clip
    ]{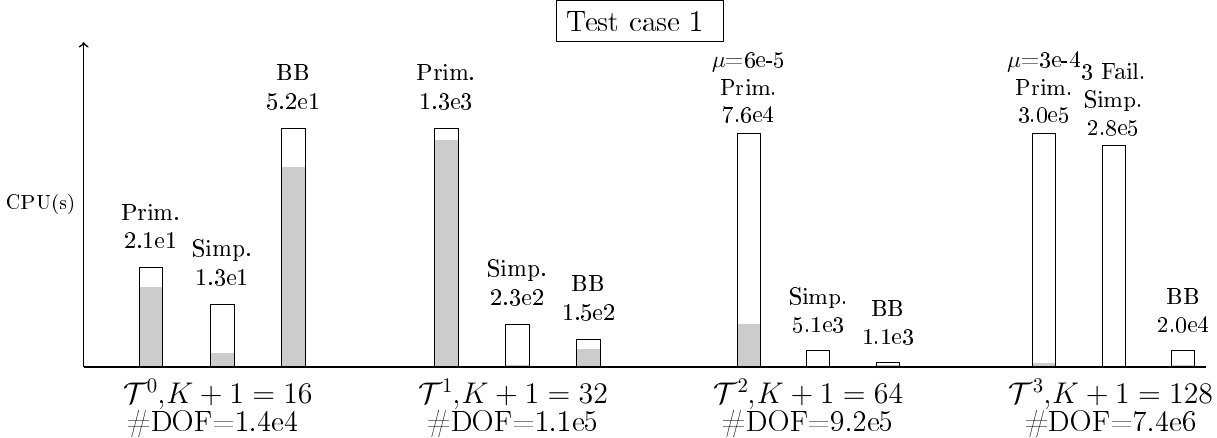}
  }
  \centerline{
    \includegraphics[width=0.8\textwidth,
      trim={0.0cm 0.0cm 0cm 0.0cm},clip
    ]{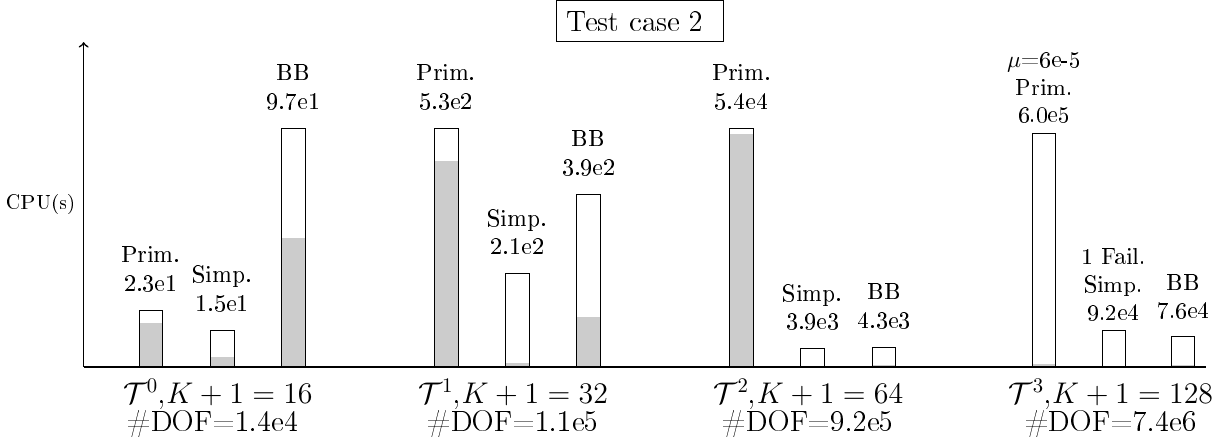}
  }
  \centerline{
    \includegraphics[width=0.8\textwidth,
    trim={0.0cm 0.0cm 0cm 0.0cm},clip
    ]{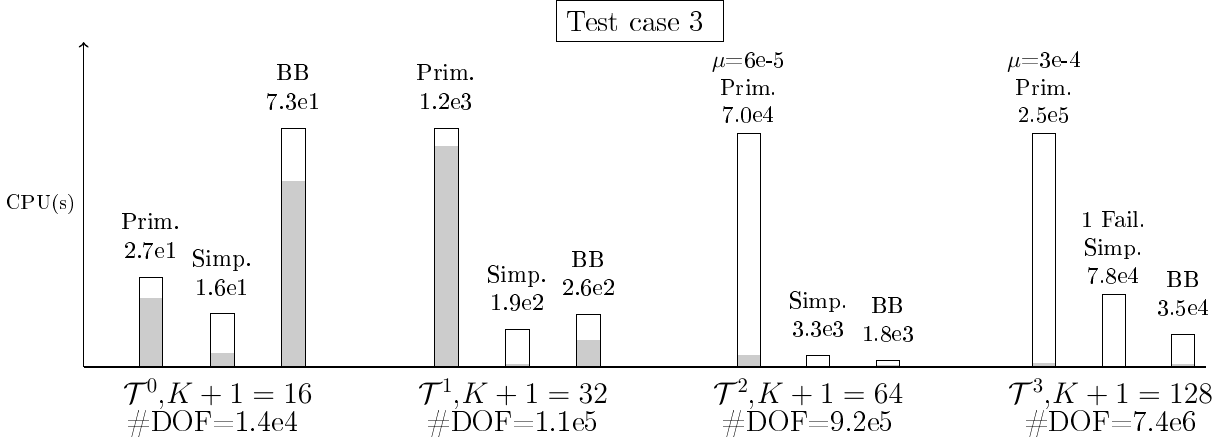}
  }
  \centerline{
    \includegraphics[width=0.65\textwidth,
    trim={0.0cm 0.0cm 0cm 0.0cm},clip
    ]{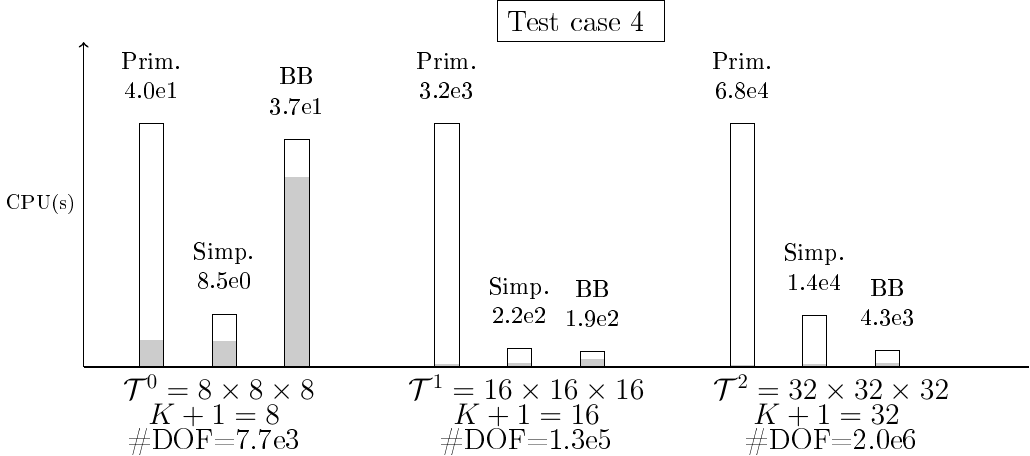}
  }
  \caption{
    Comparison of CPU time spent in solving linear algebra problems to reach an \IP{} relaxation $\Relax\approx 1e-5$ using the preconditioners approaches presented in this paper. They are denoted by Prim. (the preconditioner based on the primal Schur complement in~\cref{sec:primal}, Simp. (the SIMPLE preconditioner in~\cref{sec:dual-schur}) and BB (the $\BB$-preconditioner). The results refer to the test cases 1,2,3 and 4 (from top to bottom panel) and using different time-space discretization (left to right). The columns height is normalized by the maximum among the three preconditioners. The gray portion of the column denotes the part of the preprocessing time of each preconditioner. If $\Relax\approx 1e-5$  was not achieved, we report the value $\Relax$ reached. For the SIMPLE preconditioner some intermediate \IP{} steps failed. We report the failures number and we do not sum any CPU time to the total cost, since the $\BB$-preconditioner appears to be the best among the three presented anyway.
  }
  \label{fig:resume}
\end{figure}
In~\cref{fig:resume} we compare the total CPU time (y-axis) required to achieve \IP{} relaxation $\Relax \approx 1e-5$ (eight \IP{} iterations) with respect to the total number of degrees of freedom used, while halving the time step and the mesh size used. This summarizes the data in~\cref{tab:sol10perc,tab:sol11,tab:sol20} with gray background color. We included in this comparison the results obtain for the test cases 1, 3 and 4. 

The behavior of the preconditioner based on the primal Schur complement described in~\cref{sec:primal} is determined by the inner solver used to solve the block linear system~\ref{eq:primal-linsys}, which is the finite-dimensional counterpart of a time-space weighted and anisotropic Laplacian. This linear system was solved efficiently with the AGMG solver only for the smallest test cases during  the initial \IP{} steps, when the relaxation $\Relax$ is relatively high. In fact, it did not reach the threshold $\Relax \approx 1e-5$ in the finest discretizations for all test cases.

The SIMPLE preconditioner is rather robust with respect to different time steps and relaxation parameters $\Relax$. However, it becomes inefficient when large grids are used. In fact, the number of outer iterations and the CPU time approximately scale linearly with respect to the number of time steps and the number of interior point steps, but quadratically with respect to the number of cells. Moreover, it had some failures in the first \IP{} steps for all test cases in the finest discretization. However, the SIMPLE preconditioner is the only one that becomes faster in the last \IP{} iterations. Thus, it may be the only one providing a viable option when accurate solutions of the optimal transport problem are required, possibly combined with more efficient approaches for the initial \IP{} steps.

The $\BB$-preconditioner described in~\cref{sec:commute} is the only one able to tackle efficiently large scale problems, since the number of inner iterations required per each Krylov step remains rather constant, while the number of outer iterations increases only slightly using smaller time steps and finer grids. This holds as long as the \IP{} tolerance is approximately $1e-5$. Using smaller values, the number of inner and outer iterations tends to increase, leading in some cases to the solver's failure. The reason causing this performance degradation is not clear and needs further investigations.

\section{Conclusions}
\label{sec:conclusions}
We presented different preconditioners for solving via iterative methods the saddle point linear systems arising when computing solutions of the Benamou-Brenier formulation using an \IP{} method. The first preconditioner presented, based on the primal Schur complement, is neither efficient nor robust since its application requires the solution of a complicated linear system involving an anisotropic weighted Laplacian in the time-space domain. The second preconditioner, called SIMPLE, was shown to be remarkably robust but with a CPU time that scales quadratically with respect to the number of spatial degrees of freedom. The most efficient approach turned out to be a block triangular preconditioner where the inverse of the dual Schur complement is approximated exploiting a partial commutation of its components, which we named $\BB$-preconditioner. Although this approach loses efficiency in the latest \IP{} steps, it is the only one scaling well with respect to the time and space discretization size, both in 2d and 3d problems. Up to a reasonable value $\Relax\approx 1e-5$ for the \IP{} relaxation parameter, our numerical experiments showed that the problem can be solved with a good scaling of the CPU time with respect to the number of degrees of freedom.
Combining it with further tuning strategies of the \IP{} method (e.g. predictor-corrector methods, multilevel approach, parallelization) we believe it has the potential to provide highly efficient solvers for the dynamical optimal transport problem.
Moreover, we remark that the same strategy can be applied to variations of the problem which consider further penalization of the density in~\cref{prob:min-kinetic} (such as those considered in~\cite{Porretta2022}) without major modifications.
The deterioration of the preconditioner for smaller values of $\Relax$ requires further investigation and it is left for future works.

\section*{Acknowledgments}
The bulk of this work was completed while E.F was affiliated with Univ. Lille, Inria, CNRS, UMR 8524 - Laboratoire Paul Painlev\'e, as member of the RAPSODI team. E.F. was partially funded by the Marie Sk\l odowska-Curie Action (HE MSCA PF 101103631).   A.N. was partially funded by the Labex CEMPI (ANR-11-LABX-0007-01). The work of M. B. was supported in part by PNRR MUR Project PE0000013-FAIR.

\bibliographystyle{siam}  
\bibliography{Strings,biblio_paper}

\end{document}